\documentclass{article}
\usepackage{fancyhdr}
\usepackage[british]{babel}

\usepackage{amssymb, amsmath, amsthm, amsfonts, enumerate}
\usepackage{amsmath,setspace,scalefnt}
\usepackage[usenames,dvipsnames,svgnames,table]{xcolor}
\usepackage{graphicx,tikz,caption,subcaption}
\usetikzlibrary{patterns}
\usetikzlibrary{shapes}
\usepackage{fullpage}
\usepackage{hyperref}
\usepackage{enumitem,verbatim}
\usepackage{multicol}
\usepackage{float}
\usepackage{cleveref}
\usepackage[normalem]{ulem}
\usepackage{stackrel}

\usepackage[margin=2.5cm]{geometry}

\definecolor{gray}{rgb}{0.25, 0.25, 0.25}

\newtheorem{theorem}{Theorem}
\newtheorem{lemma}{Lemma}[section]
\newtheorem{cor}{Corollary}[section]

\newtheorem*{theorem*}{Theorem}

\theoremstyle{definition}

\newcounter{propcounter}

\theoremstyle{plain}
\newtheorem{claim}[theorem]{Claim}

\title{Berge Hamilton cycles in a random sparsification of dense hypergraphs}
\author{Seonghyuk Im\thanks{Center for AI and Natural Sciences, Korea Institute for Advanced Study (KIAS), Seoul, South Korea and Extremal Combinatorics and Probability Group (ECOPRO), Institute for Basic Science (IBS), Daejeon, South Korea. Email: {\tt seonghyuk@kias.re.kr}. SI was supported by the National Research Foundation of Korea (NRF) grant funded by the Korean government(MSIT) No. RS-2023-00210430 and by the Institute for Basic Science (IBS-R029-C4)}, Minseo Kim\thanks{Department of Mathematical Sciences, KAIST, Daejeon, South Korea. Email: {\tt minseo.kim00@kaist.ac.kr}}}
\date{\today}

\begin{document}

\maketitle

\begin{abstract}
In the standard random graph process, edges are added to an initially empty graph one by one uniformly at random. A classic result by Ajtai, Koml\'os, and Szemer\'edi, and independently by Bollob\'as, states that in the standard random graph process, with high probability, the graph becomes Hamiltonian exactly when its minimum degree becomes $2$; this is known as a \emph{hitting time} result. 
Johansson extended this result by showing the following: For a graph $G$ with $\delta(G) \geq (1/2+\varepsilon)n$, in the random graph process constrained to the host graph $G$, the hitting times for minimum degree $2$ and Hamiltonicity still coincide with high probability.

In this paper, we extend Johansson's result to Berge Hamilton cycles in hypergraphs. We prove that if an $r$-uniform hypergraph $H$ satisfies either $\delta_1(H) \geq (\frac{1}{2^{r-1}} + \varepsilon)\binom{n-1}{r-1}$ or $\delta_2(H) \geq \varepsilon n^{r-2}$, then in the random process generated by the edges of $H$, the time at which the hypergraph reaches minimum degree $2$ coincides with the time at which it contains a Berge Hamilton cycle with high probability. 
In addition, we prove an analogous result for weak Berge Hamilton cycles.
This generalizes the work of Bal, Berkowitz, Devlin, and Schacht, who established the result for the case where $H$ is a complete $r$-uniform hypergraph.

\end{abstract}

\section{Introduction}
The Erd\H{o}s-R{\'e}nyi random graph model $G(n, p)$ is a probability space of $n$-vertex graphs where each edge is chosen independently with probability $p$.
The problem of Hamiltonicity of $G(n, p)$ is one of the oldest problems in random graph theory. 
P\'{o}sa~\cite{posa}, and independently Kor\v{s}unov~\cite{korshunov}, proved that if $p = \omega(\log n/n)$, then $G(n, p)$ is Hamiltonian with high probability (w.h.p.)\footnote{We say that a sequence of events $\mathcal{E}_n$ holds with high probability if $\lim_{n \to \infty} \mathbb{P}(\mathcal{E}_n) =1$.}.
Koml\'{o}s and Szemer\'{e}di~\cite{komlos} improved this result by showing that if $p = (\log n + \log \log n + \omega(1))/n$, then $G(n, p)$ is Hamiltonian w.h.p., while if $p = (\log n + \log \log n - \omega(1))/n$, then w.h.p. $G(n, p)$ is not Hamiltonian.
Indeed, the threshold for Hamiltonicity in $G(n, p)$ coincides with the threshold for having minimum degree at least two.
Bollob\'{a}s~\cite{bollobasfirst}, and independently Ajtai, Koml\'os, and Szemer\'edi~\cite{ajtai}, proved that this minimum degree condition is an essential barrier to Hamiltonicity in the random graph by showing the `hitting time' result.

For a graph $G$ and an edge ordering $\sigma$ of $E(G)$, the \emph{subgraph process} $(G_t(\sigma))_{t=0}^{|E(G)|}$ is defined as $G_t(\sigma) = (V(G), \{\sigma(1), \sigma(2), \ldots, \sigma(t)\})$ for each $0 \leq t \leq |E(G)|$.
A \emph{random subgraph process} $(G_t)_{t=0}^{|E(G)|}$ on $G$ is a subgraph process where the edge ordering $\sigma$ is chosen uniformly at random among all $|E(G)|!$ possible orderings.
For a graph property $\mathcal{P}$, the \emph{hitting time} $\tau_{\mathcal{P}} := \tau_{\mathcal{P}}(G)$ is the smallest integer $t$ such that $G_t \in \mathcal{P}$.
When $\mathcal{P}$ is the property of having minimum degree at least $k$, we simply denote it by $\tau_k$.
Let $\mathcal{H}$ be the property of containing a Hamilton cycle.  
Then the hitting time result by Ajtai, Koml\'os, and Szemer\'edi~\cite{ajtai}, and Bollob\'as~\cite{bollobasfirst} can be stated as follows.
\begin{theorem}[\cite{ajtai, bollobasfirst}]\label{thm:hitting_complete}
    If $G = K_n$, then w.h.p. $\tau_2(G) = \tau_{\mathcal{H}}(G)$.
\end{theorem}

Afterwards, the case where $G$ is not complete was studied.
Bollob\'as and Kohayakawa~\cite{kohayakawa} proved that the hitting time statement holds when $G$ is the complete bipartite graph $K_{n, n}$. Frieze and Krivelevich~\cite{friezepseudo} and Alon and Krivelevich~\cite{alon} considered the case where $G$ is a pseudorandom graph.
Johansson~\cite{johansson} proved that the hitting time statement holds when $\delta(G) \geq (\frac{1}{2} + \varepsilon)n$ for some fixed $\varepsilon >0$.

In the context of hypergraphs, several different notions of Hamilton cycles have been studied.
For instance, Narayanan and Schacht~\cite{narayanan} determined the sharp threshold for a Hamiltonian $\ell$-cycle in an $r$-uniform hypergraph for every $r>\ell \geq 2$, and recently, Frieze and P\'{e}rez-Gim\'{e}nez~\cite{perez} confirmed the case where $\ell=1$.
See~{\cite[Section 8]{friezesurvey}} for more results.
However, the corresponding hitting time result is still unknown.

In the case of (weak) Berge Hamilton cycles, more results are known. 
For an $r$-uniform hypergraph (or an $r$-graph) $H$, a \emph{Berge cycle} is a sequence $(v_1, e_1, v_2, e_2, \ldots, v_k, e_k, v_{k+1})$ where $v_1, v_2, \ldots, v_k$ are distinct vertices with $v_1=v_{k+1}$, $e_1, e_2, \ldots, e_k$ are distinct edges, and $v_i, v_{i+1} \in e_i$ for each $1 \leq i \leq k$.
A Berge cycle is called a \emph{Berge Hamilton cycle} if $\{v_1, \ldots, v_k\} = V(H)$.
If we drop the condition that all the edges are distinct, then it is called a \emph{weak Berge cycle}.
A \emph{Berge path} and a \emph{weak Berge path} are defined similarly except that $v_1 \neq v_{k+1}$.
Poole~\cite{poole} proved a sharp threshold result for the existence of a weak Berge Hamilton cycle in a random hypergraph.
Bal, Berkowitz, Devlin, and Schacht~\cite{deepak} proved a corresponding hitting time result for both the weak Berge Hamilton cycle and the Berge Hamilton cycle. 
For an $r$-graph $H$, we define the random subgraph process $(H_t)_{t=0}^{|E(H)|}$ and the hitting time $\tau_{\mathcal{P}}(H)$ similarly to the graph case.
Let $\mathcal{BH}$ be the property of having a Berge Hamilton cycle and let $w\mathcal{BH}$ be the property of having a weak Berge Hamilton cycle.
\begin{theorem}[\cite{deepak}]\label{thm:hitting_hyper_complete}
    If $H$ is the complete $r$-graph on $n$ vertices, then w.h.p., $\tau_1(H) = \tau_{w\mathcal{BH}}(H)$ and $\tau_2(H) = \tau_{\mathcal{BH}}(H)$.
\end{theorem}

One natural direction is to consider the case where $H$ is not a complete hypergraph. 
In this paper, we extend \Cref{thm:hitting_hyper_complete} in the same spirit as Johansson's result~\cite{johansson} in the graph case.
\begin{theorem}\label{theorem: berge}
    Let $\varepsilon>0$ be a fixed constant and $r \geq 3$.
    Let $H$ be an $n$-vertex $r$-graph such that for every $v \in V(H)$, there exist at least $(\frac{1}{2} + \varepsilon)n$ vertices $u$ such that the codegree $d_H(u, v) \geq \varepsilon n^{r-2}$.
    Then w.h.p. $\tau_1(H) = \tau_{w\mathcal{BH}}(H)$ and $\tau_2(H) = \tau_{\mathcal{BH}}(H)$. 
\end{theorem}
This implies the minimum degree condition for the hitting time statement.
Let $\delta_\ell(H)$ be the \emph{minimum $\ell$-degree} of $H$, which is the minimum codegree among all $\ell$-subsets of $V(H)$.
\begin{cor}
    Let $\varepsilon>0$ be a fixed constant and $r \geq 3$.
    Let $H$ be an $n$-vertex $r$-graph satisfying one of the following two conditions.
    \begin{itemize}
        \item $\delta_1(H) \geq (\frac{1}{2^{r-1}} + \varepsilon) \binom{n-1}{r-1}$ or
        \item $\delta_2(H) \geq \varepsilon n^{r-2}$. 
    \end{itemize}
    Then w.h.p., $\tau_1(H) = \tau_{w\mathcal{BH}}(H)$ and $\tau_2(H) = \tau_{\mathcal{BH}}(H)$.
\end{cor}
\begin{proof}
    The second condition clearly implies the condition in \Cref{theorem: berge}, so we only prove that the first condition implies the condition in \Cref{theorem: berge}.
    Let $\varepsilon' = \varepsilon/(10r!)$.
    Suppose that $\delta_1(H) \geq (\frac{1}{2^{r-1}} + \varepsilon) \binom{n-1}{r-1}$.
    We choose an arbitrary vertex $v \in V(H)$ and partition $V(H) \setminus \{v\}$ into $X = \{u \in V(H) : d_H(u, v) \geq \varepsilon' n^{r-2}\}$ and $Y = V(H) \setminus (X \cup \{v\})$.
    Suppose that $|X| < (\frac{1}{2} + \varepsilon')n$.
    Then the number of edges containing $v$ and meeting $Y$ is at most $|Y| \cdot \varepsilon' n^{r-2} \leq \varepsilon' n^{r-1}$.
    All the other edges containing $v$ are contained in $X \cup \{v\}$, so the number of such edges is at most $\binom{|X|}{r-1}$.
    Thus, the degree of $v$ is at most 
        $$\binom{|X|}{r-1} + \varepsilon' n^{r-1} \leq  \binom{(\frac{1}{2} + \varepsilon')n}{r-1} + \varepsilon' n^{r-1} < (\frac{1}{2^{r-1}} + \varepsilon) \binom{n-1}{r-1},$$
    a contradiction.
    Thus, $|X| \geq (\frac{1}{2} + \varepsilon')n$. As $v$ is arbitrary, $H$ satisfies the condition in \Cref{theorem: berge}.
\end{proof}

The constant $\frac{1}{2^{r-1}}$ is best possible since if $H$ consists of two disjoint complete $r$-graphs of size $n/2$, then $H$ does not contain a Berge Hamilton cycle and has minimum degree $\binom{\lfloor(n-1)/2\rfloor}{r-1} = (\frac{1}{2^{r-1}} + o(1))\binom{n-1}{r-1}$.
Note that Kostochka, Luo, and McCourt~\cite{kostochka} proved that if $n \geq 2r+1$ and $\delta_1(H) \geq \binom{\lfloor(n-1)/2\rfloor}{r-1}+1$, then $H$ contains a Berge Hamilton cycle.
However, this minimum degree condition is not sufficient to guarantee the hitting time statement, as the following construction shows.
Assume that $n$ is divisible by $6$ and let $H$ be a $3$-graph on the vertex set $V_1 \cup V_2$ where $|V_1| = |V_2| = n/2$.
The edge set consists of complete $3$-graphs on $V_1$ and $V_2$ together with a perfect matching $M$ such that each edge of $M$ meets both $V_1$ and $V_2$ in at least one vertex.
Note that $M$ contains exactly $n/3$ edges and $\tau_2(H) = \Theta(\log n/n^2)\cdot|E(H)|$ w.h.p. (see Section~\ref{sec:value_of_tau} for details).
Thus, $H_{\tau_2}$ does not contain any edge of $M$ w.h.p., and consequently does not contain a Berge Hamilton cycle.
Therefore, the $\varepsilon$ term in the minimum degree condition is necessary for the hitting time statement.

\textbf{Organization.} In \Cref{sec:preliminaries}, we collect the definitions and the theoretical background used for the proof. In addition, we estimate the value of $\tau_2$ by applying standard computations. In \Cref{sec:proof}, we prove \Cref{theorem: berge} for the Berge Hamilton cycle. As the proof for the weak Berge Hamilton cycle is similar, we only provide a sketch of the proof for the weak Berge Hamilton cycle in \Cref{sec:weak_berge}.

\section{Preliminaries}\label{sec:preliminaries}
In this section, we collect the definitions and the theoretical background used for the proof. 
For an $r$-graph $H$ and $S \subseteq V(H)$, we denote the \emph{neighborhood of $S$} in $H$ by $N_H(S) := \{v \in V(H) \setminus S \mid \exists e \in E(H) \text{ such that } v \in e \text{ and } e \cap S \neq \emptyset\}$.
We denote the degree of $S$ in $H$ by $d_H(S)$, which is the number of edges containing $S$. 
If $S = \{v\}$ or $S = \{u, v\}$, we simply write $d_H(v)$ or $d_H(u, v)$, respectively.

We now define a booster pair, which was introduced by Montgomery~\cite{montgomery} and has been previously used for the study of hitting time by Alon and Krivelevich~\cite{alon}.
Let $H$ be an $r$-graph without a Berge Hamilton cycle. 
A pair $(e_1, e_2)$ of non-edges $e_1, e_2\in\binom{V(H)}{r}\setminus E(H)$ is called a \emph{booster pair}, or simply a \emph{BP}, if the graph $H'$ with edge set $E(H')=E(H)\cup \{e_1, e_2\}$ is either Berge Hamiltonian or has a longer Berge path than $H$ does.

\subsection{The Chernoff bound}
Let $X\sim \mathrm{Bin}(n, p)$ be a binomial random variable and $\mu=np$ be the expectation of $X$. 
Then the Chernoff bound gives an upper bound on the probability that $X$ deviates far from $\mu$. We state the standard forms of the Chernoff inequality (see, for example, \cite[Theorem 4.4]{upfal}).
\begin{theorem}
    For every $\delta>0$, the upper tail satisfies
\begin{align}\label{chernoff-U1}
    \mathbb{P}[X\geq \mu(1+\delta)]\leq \exp\left(-\mu((1+\delta)\log{(1+\delta)}-\delta)\right),
\end{align}
or simply 
\begin{align}\label{chernoff-U2}
    \mathbb{P}[X\geq \mu(1+\delta)]\leq \exp\left(-\frac{\delta^2\mu}{2+\delta}\right).    
\end{align}
For every $0<\delta<1$, the lower tail satisfies
\begin{align}\label{chernoff-L1}
\mathbb{P}[X\leq\mu(1-\delta)]\leq \exp\left(-\frac{\delta^2\mu}{2}\right).
\end{align}
\end{theorem}

\subsection{P\'osa rotation}

P\'{o}sa rotation~\cite{posa} is a powerful tool to find a Hamilton cycle in various graphs. 
We collect its hypergraph version which was introduced by Bal, Berkowitz, Devlin, and Schacht \cite{deepak}.

Let $P=(v_1, e_1, v_2, \dots, e_{\ell-1}, v_\ell)$ be a longest Berge path in a hypergraph $H$. 
Let $e$ be an edge containing $v_\ell$. 
If $e\in P$, let $e=e_i$. 
As $P$ is a Berge path, $v_i, v_{i+1}\in e$. 
Thus, $P'$ defined below is a Berge path. 
\[P'=(v_1, e_1, \dots, v_i, e, v_\ell, e_{\ell-1}, v_{\ell-1}, \dots, e_{i+1}, v_{i+1}).\] 
We now assume that $e\notin P$. If $e$ does not contain any vertex of $P$ except $v_\ell$, this contradicts the assumption that $P$ is a longest (weak) Berge path in $H$. 
Thus, $e$ must contain some vertex $v_i\in P$ with $i\neq \ell$. 
Then the $P'$ defined below is a (weak) Berge path. \[P'=(v_1, e_1, \dots, v_i, e, v_\ell, e_{\ell-1}, v_{\ell-1}, \dots, e_{i+1}, v_{i+1}).\] 
In both cases, $P'$ is a longest Berge path in $H$ with one endpoint $v_1$.
We call such operations \emph{path rotation fixing $v_1$} or simply \emph{path rotation}.

We now define an expander notion for hypergraphs and see the connection with P\'osa rotation.
An $r$-graph $H$ is a \emph{${(k, \alpha)}$-expander} if for all disjoint sets of vertices $X$ and $Y$ satisfying $|Y|<\alpha|X|$ and $|X|\leq k$, there exists an edge $e$, such that $|e\cap X|=1$ and $|e\cap Y|=0$.
The following lemma is not explicitly stated in~\cite{deepak}, but proved in the proof of Lemma 7 of~\cite{deepak}.
\begin{lemma}[{\cite{deepak}}]\label{lem21}
 Let $r \geq 3$ and let $H$ be an $r$-graph which is a $(k, 2)$-expander. 
Let $P=(v_1, e_1, v_2, \dots, e_{\ell-1}, v_\ell)$ be a longest Berge path with one end point $v_1=x$. 
Then at least $k$ distinct longest paths with distinct endpoints can be obtained from a sequence of path rotations fixing $x$.
\end{lemma}

\subsection{Asymptotic equivalence of random models}

Let $H_p$ be the space of random subgraphs of $H$ in which each edge of $H$ is included with probability $p$, independently. Let $H_t$ be the space of random subgraphs of $H$ with exactly $t$ edges from $E(H)$, chosen uniformly from the $\binom{|E(H)|}{t}$ possibilities. 
In this paper, we need to bound the probability that certain properties hold in $H_t$. But using model $H_p$ is easier for some proofs than dealing with $H_t$. The following lemma gives a relation between the probabilities in the two models.
In \cite{frieze}, they only mention the case where $H$ is a complete graph, but the proof works for a general hypergraph (see also~\cite{johansson}).

\begin{lemma}[{\cite[Lemma 1.2, 1.3]{frieze}}]\label{lem: equiv1}
Let $H$ be a hypergraph on $n$ vertices and $m$ edges. Let $\mathcal{P}$ be a graph property, $0\leq t\leq m$ and $p=t/m$. If $t, m-t\rightarrow\infty$ as $n\rightarrow\infty$,
\begin{align*}
\mathbb{P}[H_t\notin\mathcal{P}]\leq 10\sqrt{m}\cdot \mathbb{P}[H_p\notin\mathcal{P}].
\end{align*}
If $\mathcal{P}$ is a monotone increasing property, $p=o(1)$ and $mp, m(1-p)/\sqrt{mp}\rightarrow\infty$ whenever $n\rightarrow\infty$, 
\begin{align*}
\mathbb{P}[H_t\notin\mathcal{P}]\leq 3\cdot \mathbb{P}[H_p\notin\mathcal{P}].
\end{align*}
\end{lemma}

\begin{lemma}[{\cite[Corollary 1.16]{janson}}]\label{lem: equiv2}
Let $\mathcal{P}$ be a monotone increasing or decreasing property, and let $t=t(n)\rightarrow\infty$. Then $\mathbb{P}[H_p\in\mathcal{P}]\rightarrow 1$ implies $\mathbb{P}[H_t\in\mathcal{P}]\rightarrow 1$.
\end{lemma}

\subsection{Approximations of $\tau_1$ and $\tau_2$}\label{sec:value_of_tau}

To prove \Cref{theorem: berge}, we first estimate the value of  $\tau_2$ to prove the properties of $H_{\tau_2}$. 
Let $G^{(r)}(n, p)$ be a binomial random $r$-graph on $n$ vertices where each edge is included with probability $p$, independently.

\subsubsection{General bound}
We first provide a bound for $\tau_2$ independent of the structure of $H$.
Let $N=|E(H)|$, and let \[p_1=\frac{\log{n}}{n^{r-1}}, p_2=\frac{2\log{n}}{\varepsilon\cdot n^{r-1}}, m_1=N\cdot p_1, m_2=N\cdot p_2.\]

\begin{lemma}\label{lem23}
    With high probability, $m_1\leq \tau_1\leq \tau_2 \leq m_2$.
\end{lemma}

\begin{proof}
As the property that a hypergraph has minimum degree at least $k(=1$ or $2)$ is a monotone increasing property, by applying \Cref{lem: equiv1}, it is enough to prove $\mathbb{P}[\delta(H_{p_1})\geq 1]=o(1)$ and $\mathbb{P}[\delta(H_{p_2})<2]=o(1)$. For simplicity, we prove the even stronger result that $\mathbb{P}[\delta(H_{p_1})\geq 1]=o(1)$.

To prove $\mathbb{P}[\delta(H_{p_1})\geq 1]=o(1)$, we prove $\mathbb{P}[\delta(G^{(r)}(n, p_1))\geq 1]=o(1)$. 
For each vertex $v \in V(G^{(r)}(n, p_1))$, let $I_v$ be the indicator variable for the event that $v$ is isolated,  and let $X=\sum_{v}I_v$. 
The expectation of $X$ is
\begin{align*}
\mathbb{E}[X]=\sum_{v\in V(G)}\mathbb{E}[I_v]=n(1-p_1)^{\binom{n-1}{r-1}}.
\end{align*}
If $u\neq v$, $\mathbb{E}[I_v\cdot I_u]=(1-p_1)^{2\binom{n-1}{r-1}}\cdot(1-p_1)^{-\binom{n-2}{r-2}}=(1+o(1))\cdot(1-p_1)^{2\binom{n-1}{r-1}}=(1+o(1))\cdot\mathbb{E}[I_v]\cdot\mathbb{E}[I_v]$. 
Thus, the expectation of $X^2$ is
\begin{align*}
\mathbb{E}[X^2]&=\sum_{v\in V(G)}\mathbb{E}[I_v^2]+\sum_{v\neq u}\mathbb{E}[I_v\cdot I_u]=(1+o(1))\left(\sum_{v\in V(G)}\mathbb{E}[I_v^2]+\sum_{v\neq u}\mathbb{E}[I_v]\cdot\mathbb{E}[I_u]\right)\\
&=(1+o(1))\left(\sum_{v\in V(G)}\mathbb{E}[I_v]^2+\sum_{v\neq u}\mathbb{E}[I_v]\cdot\mathbb{E}[I_u]\right)=(1+o(1))\left(\sum_{v\in V(G)}\mathbb{E}[I_v]\right)^2\\
&=(1+o(1))\mathbb{E}[X]^2,
\end{align*}
where the third equality holds because $\sum_{v\in V(G)}\mathbb{E}[I_v^2]=O(n)$ and $\sum_{v\neq u}\mathbb{E}[I_v]\cdot\mathbb{E}[I_u]=\omega(n)$.
Therefore, the variance of $X$ is $o(\mathbb{E}(X)^2)$. 
By Chebyshev's inequality, $\mathbb{P}[\delta(H_{p_1})\geq 1]=o(1)$.

We now prove $\mathbb{P}[\delta(H_{p_2})<2]=o(1)$. For each vertex $v \in V(H)$, the probability that $v$ has degree less than $2$ in $H_{p_2}$ is
\begin{align*}
        \mathbb{P}[d_{H_{p_2}}(v)<2]&=(1-p_2)^{d_H(v)}+d_H(v)\cdot p_2(1-p_2)^{d_H(v)-1}\\
        & \leq (1+d_H(v)\frac{p_2}{1-p_2})\cdot(1-p_2)^{(1/2+\varepsilon)\varepsilon n^{r-1}}\\
        & \leq \frac{3}{\varepsilon}\log{n}\cdot \exp\left(-p_2(1/2+\varepsilon)\varepsilon n^{r-1}\right)\\
        & \leq \frac{3}{\varepsilon}\log{n}\cdot n^{-1-2\varepsilon}\\
        & =\frac{o(1)}{n}.
\end{align*}
By the union bound, $\mathbb{P}[\delta(H_{p_2})<2]=o(1)$ holds.
\end{proof}

\subsubsection{Tighter bound}

The bound in \Cref{lem23} is sufficient for most of the proofs in this paper. However, we need a tighter bound for some parts of the proof.
We note that the calculations in this section were inspired by \cite{johansson}. 
Let $p_0$ be the unique real root of the following equation:
\begin{align*}
\sum_{v\in V(H)}(1-p)^{d_H(v)}=(\log{n})^{-1}.
\end{align*}
Note that there exists a unique real root as the left hand side is strictly decreasing as $p$ increases. Let $\gamma = \gamma(n)$ be a nonnegative increasing function such that $\gamma = \Theta(\log{\log{\log{n}}})$.
Let 
\begin{align*}
p_3=p_0-\frac{\gamma}{n^{r-1}}, p_4=p_0+\frac{\gamma}{n^{r-1}}, m_3=N\cdot p_3, m_4=N\cdot p_4.
\end{align*}
We note that if we substitute $p=0.5(r-1)!\log{n}/n^{r-1}$, then $\sum_{v\in V(H)}(1-p)^{d_H(v)}\geq n(1-p)^{n^{r-1}/(r-1)!}>1$ and if we substitute $p=3\log{n}/\varepsilon n^{r-1}$, then $\sum_{v\in V(H)}(1-p)^{d_H(v)}\leq n \cdot \exp(-0.5p\varepsilon n^{r-1})=n^{-0.5}$. Thus, $p_0=\Theta(\log{n}/n^{r-1})$. Moreover, we can deduce
\begin{align}\label{eqn: p_3}
p_3\geq\frac{\log{n}}{4n^{r-1}},
\end{align}
since $r\geq 3$ and
\begin{align}\label{eqn: p_4}
p_4\leq\frac{4\log{n}}{\varepsilon n^{r-1}},
\end{align}
provided that $n$ is large enough.
To prove that $m_3\leq \tau_2\leq m_4$ w.h.p., we need the following lemma.
\begin{lemma}\label{lem24}
The following inequalities hold:
\begin{enumerate}[label=\upshape \textbf{(L\arabic{enumi})}]
    \item\label{L1}
    $\log{n}\cdot\sum_{v\in V(H)}(1-p_3)^{d_H(v)}\geq e^{\varepsilon\gamma/4}$.
    \item\label{L2}
    $\log{n}\cdot\sum_{v\in V(H)}(1-p_3)^{d_H(v)}\leq e^{\gamma}$.
    \item\label{L3}
    $\log{n}\cdot\sum_{v\in V(H)}(1-p_4)^{d_H(v)}\leq e^{-\gamma}$.
\end{enumerate}
\end{lemma}

\begin{proof}
\begin{description}
\item{\ref{L1}}
For each $v\in V(H)$, we have  
\begin{align*}
(1-p_3)^{d_H(v)}&=(1-p_0)^{d_H(v)}\left(1+\frac{p_0-p_3}{1-p_0}\right)^{d_H(v)}\\
&\geq(1-p_0)^{d_H(v)}\left(1+\frac{\gamma}{n^{r-1}}\right)^{\varepsilon n^{r-1}/2}\\
&\geq (1-p_0)^{d_H(v)}\cdot\exp\left(\frac{\varepsilon\gamma}{4}\right).
\end{align*}
By summing over all $v\in V(H)$, we conclude
\begin{align*}
\sum_{v\in V(H)}(1-p_3)^{d_H(v)}\geq\sum_{v\in V(H)}(1-p_0)^{d_H(v)}\cdot e^{\varepsilon\gamma/4}\geq\frac{e^{\varepsilon\gamma/4}}{\log{n}}.
\end{align*}

\item{\ref{L2}}
For each $v\in V(H)$, we have  
\begin{align*}
(1-p_3)^{d_H(v)}&=(1-p_0)^{d_H(v)}\left(1+\frac{p_0-p_3}{1-p_0}\right)^{d_H(v)}\\
&\leq(1-p_0)^{d_H(v)} \cdot\exp\left(\frac{\gamma}{n^{r-1}} d_H(v)\right)\\
&\leq (1-p_0)^{d_H(v)}\cdot\exp(\gamma).
\end{align*}
By summing over all $v\in V(H)$, we conclude
\begin{align*}
\sum_{v\in V(H)}(1-p_3)^{d_H(v)}\leq\sum_{v\in V(H)}(1-p_0)^{d_H(v)}\cdot e^{\gamma}\leq\frac{e^{\gamma}}{\log{n}}.
\end{align*}

\item{\ref{L3}} It can be obtained by following the same calculation as in \ref{L2}.

\end{description}
\end{proof}

Note that \Cref{lem24} holds for any suitable choice of $\gamma$.
We now prove that $m_3$ and $m_4$ are the desired bounds for $\tau_2$.
\begin{lemma}\label{lem26}
With high probability, $m_3\leq \tau_2\leq m_4$.
\end{lemma}

\begin{proof}
By \Cref{lem: equiv1}, it is enough to prove $\mathbb{P}[\delta(H_{p_3})\geq 2]=o(1)$ and $\mathbb{P}[\delta(H_{p_4})<2]=o(1)$. 
For any vertex $v\in V(H)$, the probability that $v$ has degree less than $2$ in $H_{p_4}$ is
\begin{align*}
\mathbb{P}[d_{H_{p_4}}(v)< 2]&=(1-p_4)^{d_{H}(v)}+p_4\cdot d_H(v)(1-p_4)^{d_H(v)-1}\\
&\leq 2p_4\cdot d_H(v)(1-p_4)^{d_H(v)}\\
&\leq (8/\varepsilon)(1-p_4)^{d_H(v)}\cdot\log{n},
\end{align*}
where the first inequality holds since $\frac{1}{1-p_4}\leq 1.5$ and $p_4\cdot d_H(v)\geq 2$ for sufficiently large $n$. By the union bound and \ref{L3}, we have
\begin{align*}
\mathbb{P}[\delta(H_{p_4})<2]\leq\sum_{v\in V(H)}\mathbb{P}[d_{H_{p_4}}(v)<2]\leq(8/\varepsilon)\log{n}\cdot\sum_{v\in V(H)}(1-p_4)^{d_H(v)}=o(1).
\end{align*}

We now prove $\mathbb{P}[\delta(H_{p_3})\geq 2]=o(1)$. Let $I_v$ be an indicator variable for the event $\{d_{H_{p_3}}(v)\leq 1\}$, and let $X=\sum_{v\in V(H)}I_v$. 
Then the expectation of $X$ is at least
\begin{align*}
\mathbb{E}[X]&=\sum_{v\in V(H)}\left[(1-p_3)^{d_H(v)}+p_3\cdot d_H(v)\cdot(1-p_3)^{d_H(v)-1}\right]\\
&\geq\sum_{v\in V(H)}\left(1+\frac{p_3\cdot d_H(v)}{1-p_3}\right)\cdot(1-p_3)^{d_H(v)}\\
&\geq\Omega(\log{n})\cdot\sum_{v\in V(H)}(1-p_3)^{d_H(v)}\\
&\geq\Omega(\log{n})\cdot\frac{e^{\varepsilon\gamma/4}}{\log{n}}\\
&=\omega(1),
\end{align*}
where the third inequality holds by \ref{L1}. 

Now we investigate the expectation of $X^2$. If $u\neq v$, the expectation of $I_u\cdot I_v$ is
\begin{align*}
\mathbb{E}[I_u\cdot I_v]&=(1-p_3)^{d_H(u, v)}\left[(1-p_3)^{d_H(u)-d_H(u, v)}+p_3(d_H(u)-d_H(u, v))(1-p_3)^{d_H(u)-d_H(u, v)-1}\right]\\
&\quad\cdot\left[(1-p_3)^{d_H(v)-d_H(u, v)}+p_3(d_H(v)-d_H(u, v))(1-p_3)^{d_H(v)-d_H(u, v)-1}\right]\\
&\quad+p_3\cdot d_H(u, v)(1-p_3)^{d_H(u)+d_H(v)-d_H(u, v)-1}\\
&=(1-p_3)^{d_H(u)+d_H(v)-d_H(u, v)}\cdot\left(1+\frac{p_3}{1-p_3}(d_H(u)-d_H(u, v))\right)\\
&\quad\cdot\left(1+\frac{p_3}{1-p_3}(d_H(v)-d_H(u, v))\right)+o\left((1-p_3)^{d_H(u)+d_H(v)-d_H(u, v)}\right)\\
&=(1+o(1))(1-p_3)^{d_H(u)+d_H(v)}\cdot\frac{p_3}{1-p_3}(d_H(u)-d_H(u, v))\cdot\frac{p_3}{1-p_3}(d_H(v)-d_H(u, v))\\
&=(1+o(1))(1-p_3)^{d_H(u)+d_H(v)}\cdot\frac{p_3}{1-p_3}d_H(u)\cdot\frac{p_3}{1-p_3}d_H(v)\\
&=(1+o(1))\mathbb{E}[I_u]\mathbb{E}[I_v], 
\end{align*}
where the second and the fourth equalities come from $d_H(u, v)=o(d_H(u))$, $d_H(u, v)=o(d_H(v))$ and the third equality holds as $d_H(u, v)=O(n^{r-2})$ and $p_3=\Theta(\log{n}/n^{r-1})$. 
Moreover, since $\mathbb{E}[X]\rightarrow\infty$ as $n \to \infty$, we have $\mathbb{E}[X]=\sum_{v\in V(H)}\mathbb{E}[I_v^2] = o(\mathbb{E}[X^2])$. Therefore, expectation of $X^2$ is equal to
\begin{align*}
\mathbb{E}[X^2]&=\sum_{v\in V(H)}\mathbb{E}[I_v^2]+\sum_{v\neq u}\mathbb{E}[I_v\cdot I_u]=(1+o(1))\left(\sum_{v\in V(H)}\mathbb{E}[I_v^2]+\sum_{v\neq u}\mathbb{E}[I_v]\cdot\mathbb{E}[I_u]\right)\\
&=(1+o(1))\left(\sum_{v\in V(H)}\mathbb{E}[I_v]^2+\sum_{v\neq u}\mathbb{E}[I_v]\cdot\mathbb{E}[I_u]\right)=(1+o(1))\left(\sum_{v\in V(H)}\mathbb{E}[I_v]\right)^2\\
&=(1+o(1))\mathbb{E}[X]^2.
\end{align*}

It implies $\mathrm{Var}(X)=o(\mathbb{E}[X]^2)$. Chebyshev's inequality shows that
\begin{align*}
\mathbb{P}[X=0]\leq\mathbb{P}[|X-\mathbb{E}[X]|\geq\mathbb{E}[X]]\leq\frac{\mathrm{Var}(X)}{\mathbb{E}[X]^2}=o(1).
\end{align*}
It proves $\mathbb{P}[\delta(H_{p_3})\geq 2]=o(1)$.
\end{proof}

\section{Proof of \Cref{theorem: berge} for Berge Hamiltonicity}\label{sec:proof}
We may assume that $\varepsilon>0$ is a sufficiently small constant throughout this section.
If a hypergraph $H$ has a Berge Hamilton cycle, $H$ has a minimum degree of at least two. Thus, it is sufficient to prove $H_{\tau_{2}}$ is Berge Hamiltonian w.h.p. 
The following two lemmas are the key ingredients for the proof.

\begin{lemma}\label{lem31}
With high probability, $H_{\tau_{2}}$ contains a subgraph $\Gamma$ which is a connected $(\frac{\varepsilon n}{6}, 2)$-expander with at most $\varepsilon^8 n\log{n}+\frac{6}{\varepsilon}$ edges. 
\end{lemma}

\begin{lemma}\label{lem32}
With high probability, for every subgraph $\Gamma$ of $H_{\tau_2}$ which is a connected, non-Hamiltonian $(\frac{\varepsilon n}{6}, 2)$-expander with at most $2\varepsilon^8 n\log{n}$ edges, $H_{\tau_{2}}$ contains a booster pair with respect to $\Gamma$. 
\end{lemma}

These two lemmas immediately imply \Cref{theorem: berge} as follows.
\begin{proof}[Proof of \Cref{theorem: berge} for Berge Hamiltonicity]
By \Cref{lem31}, there exists a connected $(\frac{\varepsilon n}{6}, 2)$-expander $\Gamma_1$ in $H_{\tau_2}$ with high probability.
If $\Gamma_1$ is Hamiltonian, then we are done.
If $\Gamma_1$ is not Hamiltonian, by \Cref{lem32}, $H_{\tau_2}$ contains a booster pair with respect to $\Gamma_1$, and thus
we can add a booster pair to $\Gamma_1$ to make the longest Berge path longer. Let $\Gamma_2$ be the hypergraph obtained by adding such a booster pair to $\Gamma_1$.
Then as $\Gamma_2$ is a connected $(\frac{\varepsilon n}{6}, 2)$-expander with at most $2\varepsilon^8 n\log{n}$ edges, either $\Gamma_2$ has a Berge Hamilton cycle or $H_{\tau_2}$ contains a booster pair with respect to $\Gamma_2$.
By repeating this process at most $n$ times, we obtain a subgraph $\Gamma_k$ for some $k \in [n]$ of $H_{\tau_2}$ which contains a Berge Hamilton cycle. Therefore, $H_{\tau_2}$ is Berge Hamiltonian with high probability.
\end{proof}
In the rest of this section, we prove \Cref{lem31} and \Cref{lem32}.

\subsection{Proof of \Cref{lem31}}

To prove \Cref{lem31}, we first prove some good properties of $H_{\tau_2}$ (\Cref{lem: properties}), and then prove there exists an expander using such good properties (\Cref{lem: expander}). 
For a hypergraph $\Gamma$ on $n$ vertices, let \[d_0:=\varepsilon^8\log{n},\] and let \[\mathrm{SMALL}(\Gamma):=\{v\in V(\Gamma)\,:\,d_{\Gamma}(v)\leq d_0\}.\] 

\begin{lemma}\label{lem: properties}
The hypergraph $G=H_{\tau_2}$ has the following properties w.h.p.:
\begin{enumerate}[label=\upshape \textbf{(P\arabic{enumi})}]
    \item \label{(P1)}
$\Delta(G)\leq \frac{10}{\varepsilon}\log{n}$.
    \item \label{(P2)}
$|\mathrm{SMALL}(G)|\leq n^{0.1}$.
    \item \label{(P3)}
No edge intersects $\mathrm{SMALL}(G)$ at more than one vertex, and no $v\in V(G)\setminus \mathrm{SMALL}(G)$ lies in more than one edge incident to $N[\mathrm{SMALL}(G)]-\{v\}$.  
    \item \label{(P4)}
For any $U\subseteq V(G)$ with $|U|\leq\frac{n}{\sqrt{\log{n}}}$, there are at most $|U|(\log{n})^{3/4}$ edges of $G$ that intersect $U$ at more than one vertex.
    \item \label{(P5)}
For any disjoint vertex sets $U, W\subseteq V(G)$ satisfying $|U|\leq\frac{n}{\sqrt{\log{n}}}$ and $|W|\leq3|U|$, there are at most $0.5\varepsilon\log{n}|U|$ edges of $G$ that intersect $U$ at exactly one vertex and also intersect $W$. 
    \item \label{(P6)} 
For any disjoint vertex sets $U, W\subseteq V(G)$ satisfying $|U|=\frac{n}{\sqrt{\log{n}}}$ and $|W|=\left(1-\frac{\varepsilon}{2}\right)n$, there are at least $n(\log{n})^{1/3}$ edges of $G$ that intersect $U$ at exactly one vertex and intersect $W$ at exactly $r-1$ vertices. 
    \item \label{(P7)}
Let $U, W\subseteq V(G)$ be disjoint vertex sets whose union is $V(G)$ and which satisfy $\frac{\varepsilon n}{6}\leq|U|\leq|W|\leq n-\frac{\varepsilon n}{6}$. Then there exists at least one edge that intersects both $U$ and $W$.
\end{enumerate}
\end{lemma}

\begin{proof}
For properties which are monotone increasing or decreasing, we can use \Cref{lem: equiv2}. By \Cref{lem23}, if a property is increasing, it is enough to show that $H_{p_1}$ or $H_{p_3}$ has such a property w.h.p., and if a property is decreasing, it is enough to show that $H_{p_2}$ or $H_{p_4}$ has such a property w.h.p. 

\begin{description}
\item{\ref{(P1)}} 
We prove $\mathbb{P}[H_{p_2}\notin \textbf{(P1)}]=o(1)$. By the union bound, it is enough to show $\mathbb{P}[d_{H_{p_2}}(v)\geq(10/\varepsilon)\log{n}]=\frac{o(1)}{n}$ for every $v \in V(H)$. Let $M=n^{r-1}/(r-1)!$ and note that $d_H(v)\leq M$ for every $v \in V(H)$. By the Chernoff bound, we have
\begin{align*}
\mathbb{P}[d_{H_{p_2}}(v)\geq \frac{10}{\varepsilon}\log{n}] \leq \mathbb{P}[\mathrm{Bin}(M, p_2)\geq 5M\cdot p_2]\stackrel{\eqref{chernoff-U2}}{\leq} \exp\left(\frac{-16\log{n}}{3\varepsilon}\right)=\frac{o(1)}{n},
\end{align*}
as $n$ is large enough.
    \item{\ref{(P2)}}
To prove $\mathbb{P}[H_{p_3}\notin\textbf{(P2)}]=o(1)$, we first bound $\mathbb{P}[v\in \mathrm{SMALL}(H_{p_3})]$ for each $v\in V(H)$.
\begin{align*}
\mathbb{P}[d_{H_{p_3}}(v)\leq \varepsilon^8\log{n}]=\sum_{i=0}^{\varepsilon^8\log{n}}\binom{d_H(v)}{i}p_3^i(1-p_3)^{d_H(v)-i}.
\end{align*}
Since $\binom{d_H(v)}{i+1}p_3^{i+1}(1-p_3)^{d_H(v)-i-1}/\binom{d_H(v)}{i}p_3^i(1-p_3)^{d_H(v)-i}=\frac{p_3}{1-p_3}\cdot\frac{d_H(v)-i}{i+1}>2$ for $i<\varepsilon^8\log{n}$, by the sum of a geometric series, we have
\begin{align*}
\mathbb{P}[d_{H_{p_3}}(v)\leq \varepsilon^8\log{n}] & \leq 2\binom{d_H(v)}{\varepsilon^8\log{n}}p_3^{\varepsilon^8\log{n}}(1-p_3)^{d_H(v)-\varepsilon^8\log{n}}\\
& \leq(2+o(1))\left(\frac{e\cdot p_3\cdot d_H(v)}{\varepsilon^8\log{n}}\right)^{\varepsilon^8\log{n}}(1-p_3)^{d_H(v)}\\
& \stackrel{\eqref{eqn: p_4}}{\leq}(2+o(1))\left(\frac{4e}{\varepsilon^9}\right)^{\varepsilon^8\log{n}}(1-p_3)^{d_H(v)}\\
& =(2+o(1))e^{\varepsilon^8\log{(4e/\varepsilon^9)}}(1-p_3)^{d_H(v)}\\
&\leq 2.5e^{0.01}(1-p_3)^{d_H(v)}
\end{align*}
where the last inequality holds from the fact that $\varepsilon^8\log{(4e/\varepsilon^9)}<0.01$ when $\varepsilon$ is small enough. Thus, by \Cref{lem24}, the expectation of $|\mathrm{SMALL}(H_{p_3})|$ is
\begin{align*}
\mathbb{E}[|\mathrm{SMALL}(H_{p_3})|]\leq 2.5e^{0.01}\sum_{v\in V(H)}(1-p_3)^{d_H(v)}\leq 3n^{0.01}.
\end{align*}
By the Markov inequality, $\mathbb{P}[|\mathrm{SMALL}(H_{p_3})|\geq n^{0.1}]=o(1)$.
    \item{\ref{(P3)}} 
Since property \ref{(P3)} is neither an increasing nor decreasing property, we cannot use \Cref{lem: equiv2} directly. For hypergraphs $J_1 \subseteq J_2$, we say $J_2$ satisfies the property $\mathbf{(P3-J_1)}$ if no edge of $J_2$ intersects $\mathrm{SMALL}(J_1)$ more than once, and no $v\in V(J_2)\setminus \mathrm{SMALL}(J_1)$ lies in more than one edge of $J_2$ incident to $N[\mathrm{SMALL}(J_1)]-\{v\}$.

We claim that with high probability, $H_{m_4}$ satisfies $\mathbf{(P3-H_{m_3})}$. As the existence of a path is an increasing property and the property of a vertex being contained in $\mathrm{SMALL}$ is a decreasing property, together with the fact that $m_3 \leq \tau_2 \leq m_4$ w.h.p., this implies $H_{\tau_2}$ satisfies \ref{(P3)} with high probability.
To prove the claim, we first show that $\mathbf{(P3-H_{m_3})}$ holds for $H_{m_3}$, and the second step is to prove that after adding $m_4-m_3$ edges uniformly at random to $H_{m_3}$, $\mathbf{(P3-H_{m_3})}$ still holds w.h.p.

We now prove that $H_{m_3}$ satisfies \ref{(P3)}. For $H_{m_3}$, not satisfying $\mathbf{(P3-H_{m_3})}$ is equivalent to:
there exist $u\notin \mathrm{SMALL}(H_{m_3})$, $s, t\in \mathrm{SMALL}(H_{m_3})$, $a, b\in V(H)$, and $e_i\in E(H_{m_3})$ for $i \in [4]$, satisfying one of the following: 
\begin{enumerate}
    \item{}
    $s, u, e_1, e_2$ with $s, u\in e_1\cap e_2$;
    \item{}
    $s, t, e_1$ with $s, t\in e_1$;
    \item{}
    $s, t, u, e_1, e_2$ with $s, u\in e_1$ and $t, u\in e_2$;    
    \item{}
    $s, a, u, e_1, e_2$ with $s, a, u\in e_1$ and $a, u\in e_2$;
    \item{}
    $s, a, u, e_1, e_2, e_3$ with $s, u\in e_1$, $u, a\in e_2$, and $s, a\in e_3$;
    \item{}
    $s, t, a, u, e_1, e_2, e_3$ with $s, u\in e_1$, $u, a\in e_2$, and $t, a\in e_3$;
    \item{}
    $s, a, u, e_1, e_2, e_3$ with $s, a\in e_1$ and $u, a\in e_2\cap e_3$;
    \item{}
    $s, a, b, u, e_1, e_2, e_3$ with $s, a, b\in e_1$, $a, u\in e_2$, and $b, u\in e_3$;
    \item{}
    $s, a, b, u, e_1, e_2, e_3, e_4$ with $s, a\in e_1$, $a, u\in e_2$, $b, u\in e_3$, and $s, b\in e_4$;
    \item{}
    $s, t, a, b, u, e_1, e_2, e_3, e_4$ with $s, a\in e_1$, $a, u\in e_2$, $u, b\in e_3$, and $b, t\in e_4$;
\end{enumerate}
In each case, the listed vertices and edges are assumed to be distinct.

We first fix a tuple $T=(s, t, u, a, b, e_1, e_2, e_3, e_4) \in V(H)^5 \times E(H)^4$ that satisfies one of the above conditions.
For each $i \in [10]$, if $T$ satisfies condition $i$, let $T_i$ be a subtuple of $T$ consisting of the elements that are used in condition $i$. 
For instance, $T_1 = (s, u, e_1, e_2)$ and $T_{10}=T$. 
If $T$ does not satisfy the condition $i$, we do not define $T_i$.
Let $V(T_i)$ be the set of vertices contained in $T_i$, and let $E(T_i)$ be the set of edges contained in $T_i$.
We aim to bound the probability $\mathbb{P}[s \in \mathrm{SMALL}(H_{m_3}), E(T_i) \subseteq E(H_{m_3})]$ for each $i\in[10]$. 
We separate the proof into two cases depending on whether $T_i$ contains one or two vertices in $\mathrm{SMALL}(H_{m_3})$.
\newline

\noindent{\bf 1. Cases} $\mathbf{i\in\{1, 4, 5, 7, 8, 9\}}$.

 Let $H^1$ be a spanning subgraph of $H$ with $E(H^1)=E(H)\setminus E(T_i)$, and let $p'_3:=\frac{m_3-|E(T_i)|}{N}$. 
 We observe that the conditional event $\{s \in \mathrm{SMALL}(H_{m_3})\,|\, E(T_i) \subseteq E(H_{m_3})\}$ can be viewed as the event that the degree of $s$ is at most $d_0-|E(T_i) \cap e(s)|$ in the random hypergraph $H^1_{m_3 - |E(T_i)|}$, where $e(s)$ is the set of edges in $H$ that are incident to $s$. Thus, by \Cref{lem: equiv1}, the probability $\mathbb{P}[s \in \mathrm{SMALL}(H_{m_3}), E(T_i) \subseteq E(H_{m_3})]$ is at most
 
\begin{align*}
\mathbb{P}[s \in \mathrm{SMALL}(H_{m_3}), E(T_i) \subseteq E(H_{m_3})]&=\mathbb{P}[s \in \mathrm{SMALL}(H_{m_3})\,|\, E(T_i) \subseteq E(H_{m_3})]\\
&\cdot\mathbb{P}[E(T_i) \subseteq E(H_{m_3})]\\
&= \mathbb{P}[s \in \mathrm{SMALL}(H_{m_3})\,|\, E(T_i) \subseteq E(H_{m_3})]\cdot\frac{\binom{N-|E(T_i)|}{m_3-|E(T_i)|}}{\binom{N}{m_3}}\\
&\leq3\cdot\mathbb{P}[s\in \mathrm{SMALL}(H^1_{p'_3})]\cdot \frac{\binom{N-|E(T_i)|}{m_3-|E(T_i)|}}{\binom{N}{m_3}}\\
&\leq 3\cdot\mathbb{P}[\mathrm{Bin}\left(\frac{\varepsilon}{3}n^{r-1}, p'_3\right)\leq\varepsilon^8\log{n}]\cdot\left(\frac{m_3}{N}\right)^{|E(T_i)|}\\
&\stackrel{\eqref{chernoff-L1}}{\leq} 3\cdot\exp\left(-(0.99)^2\varepsilon\log{n}/6\right)\cdot p_3^{|E(T_i)|}\\
&\leq 3\cdot n^{-\varepsilon/12}\cdot p_3^{|E(T_i)|},
\end{align*}
where the second inequality holds since $\deg_{H^1}(s) \geq \deg_H(s)-|E(T_i)|\geq\frac{\varepsilon}{2}n^{r-1}-4\geq\frac{\varepsilon}{3}n^{r-1}$ and the third inequality comes from $p_3'\geq\frac{\log{n}}{n^{r-1}}$.

The number of ways to choose $T_i$ in $H$ is at most $n^{|E(T_i)|(r-1)}$ for each $i\in\{1, 4, 5, 7, 8, 9\}$. Since $p_3=O\left(\frac{\log{n}}{n^{r-1}}\right)$, the probability that case $i$ occurs is at most $3n^{-\varepsilon/12}\cdot O\left((\log{n})^{|E(T_i)|}\right)=o(1)$. \newline

\noindent{\bf 2. Cases} $\mathbf{i\in\{2, 3, 6, 10\}}$.

For a subset of edges $F\subseteq E(H)$, let $\Sigma(F, k)$ be the event that $|F\cap E(H_{m_3})|= k$.
Similarly, let $\Sigma(F, \geq k)$ be the event that $|F\cap E(H_{m_3})|\geq k$ and let $\Sigma(F, \leq k)$ be the event that $|F\cap E(H_{m_3})|\leq k$. For a pair of distinct vertices $v_1, v_2\in V(H)$, let $e(v_1, v_2)=e_H(v_1, v_2)$ be the set of edges in $H$ that are incident to both $v_1$ and $v_2$. The probability of $\{s, t \in \mathrm{SMALL}(H_{m_3}), E(T_i) \subseteq E(H_{m_3})\}$ is 
\begin{align*}
\mathbb{P}[s, t \in \mathrm{SMALL}(H_{m_3}), E(T_i) \subseteq E(H_{m_3})]&=\sum_{k=0}^{m_3}\mathbb{P}[s, t \in \mathrm{SMALL}(H_{m_3})\,|\,E(T_i) \subseteq E(H_{m_3}), \Sigma(e(s, t), k)]\\
&\quad\cdot\mathbb{P}[E(T_i) \subseteq E(H_{m_3}), \Sigma(e(s, t), k)]\\
&=\sum_{k=0}^{m_3}\mathbb{P}[s\in \mathrm{SMALL}(H_{m_3})\,|\,E(T_i) \subseteq E(H_{m_3}), \Sigma(e(s, t), k)]\\
&\quad\cdot\mathbb{P}[t\in \mathrm{SMALL}(H_{m_3})\,|\,E(T_i) \subseteq E(H_{m_3}), \Sigma(e(s, t), k)]\\
&\quad\cdot\mathbb{P}[E(T_i) \subseteq E(H_{m_3}), \Sigma(e(s, t), k)],
\end{align*}
where the last equality holds since the events $\{s\in \mathrm{SMALL}(H_{m_3})\}$ and $\{t\in \mathrm{SMALL}(H_{m_3})\}$ are independent conditioned on the event $\{E(T_i) \subseteq E(H_{m_3}), \Sigma(e(s, t), k)\}$.

We now prove that the codegree of $s$ and $t$ in $H_{m_3}$ is at most $10r$ with high probability. As $|e(s, t) \cap E(T_i)|, |E(T_i) \setminus e(s, t)| \leq |E(T_i)| \leq 4$, we have 

\begin{align*}
\mathbb{P}[\Sigma(e(s, t), \geq 10r+1)\,|\,E(T_i)\subseteq E(H_{m_3})]&\leq\sum_{k=10r+1}^{m_3}\frac{\binom{|e(s, t)|-|e(s, t)\cap E(T_i)|}{k-|e(s, t)\cap E(T_i)|}\binom{N-|e(s, t)|-|E(T_i)\setminus e(s, t)|}{m_3-k-|E(T_i)\setminus e(s, t)|}}{\binom{N-|E(T_i)|}{m_3-|E(T_i)|}}\\
&\leq\sum_{k=10r+1}^{m_3}\binom{|e(s, t)|}{k}\frac{\binom{N-k}{m_3-k}}{\binom{N-|E(T_i)|}{m_3-|E(T_i)|}}\\
&\leq\sum_{k=10r+1}^{m_3}\binom{|e(s, t)|}{k}\left(\frac{m_3}{N}\right)^{k-|E(T_i)|}\\
&=o\left(n^{(r-1)|E(T_i)|}\right)\sum_{k=10r+1}^{m_3}\left(\frac{e\cdot p_3\cdot|e(s, t)|}{k}\right)^k\\
&\leq o\left(n^{4(r-1)}\right)\sum_{k=10r+1}^{m_3}\left(\frac{2e\log{n}}{\varepsilon kn}\right)^{k}\\
&\leq o\left(n^{4(r-1)}\right)\cdot\left(\frac{2e\log{n}}{10\varepsilon rn}\right)^{10r}\\
&=o(n^{-10})
\end{align*}
where the sixth inequality holds since if we define $f(k)=\left(\frac{2e\log{n}}{\varepsilon kn}\right)^k$, then $f(k)/f(k+1)>2$ for $k\leq m_3$, and we use the sum of a geometric series. This result yields 
\begin{align}\label{eqn6}
\mathbb{P}[s, t \in \mathrm{SMALL}(H_{m_3}), E(T_i) \subseteq E(H_{m_3})]&\leq \sum_{k=0}^{10r}\mathbb{P}[s\in \mathrm{SMALL}(H_{m_3})\,|\,E(T_i) \subseteq E(H_{m_3}), \Sigma(e(s, t), k)] \notag \\
&\quad\cdot\mathbb{P}[t\in \mathrm{SMALL}(H_{m_3})\,|\,E(T_i) \subseteq E(H_{m_3}), \Sigma(e(s, t), k)] \notag \\
&\quad\cdot\mathbb{P}[E(T_i) \subseteq E(H_{m_3}), \Sigma(e(s, t), k)] \notag \\
&\quad+\mathbb{P}[E(T_i)\subseteq E(H_{m_3})]\cdot o(n^{-10}).
\end{align}

Since $\{s\in \mathrm{SMALL}(H_{m_3})\,|\,E(T_i) \subseteq E(H_{m_3})\}$ and $\{t\in \mathrm{SMALL}(H_{m_3})\,|\,E(T_i) \subseteq E(H_{m_3})\}$ are decreasing properties, we can use \cref{lem: equiv1} to bound the probability. We fix $0 \leq k \leq 10r$ and let $p=(m_3-k-|E(T_i)\setminus e(s, t)|)/N$. Then we have
\begin{align*}
&\mathbb{P}[s\in \mathrm{SMALL}(H_{m_3})\,|\,E(T_i) \subseteq E(H_{m_3}), \Sigma(e(s, t), k)]\\
&\leq3\cdot\sum_{a=0}^{\varepsilon^8\log{n}}\binom{d_H(s)-|e(s, t)|}{a}p^a(1-p)^{d_H(s)-|e(s, t)|-a}\\
&\leq(3+o(1))\cdot\left(1+\sum_{a=1}^{\varepsilon^8\log{n}}\left(\frac{e\cdot p\cdot d_H(s)}{a}\right)^{a}\right)\cdot(1-p)^{d_H(s)}\tag{$*$}\label{star}
\end{align*}
as $|e(s, t)|+\varepsilon^8\log{n}=O(n^{r-2})$.
Let $f(i)=\left(\frac{e\cdot p\cdot d_H(s)}{i}\right)^{i}$. If $i\leq\varepsilon^8\log{n}$, then $\frac{f(i+1)}{f(i)}>2$ whenever $\varepsilon$ is small enough. Thus, the sum of a geometric series yields:
\begin{align}\label{eqn7}
\eqref{star}&\leq\left(4+(6+o(1))\left(\frac{e\cdot p\cdot d_H(s)}{\varepsilon^8\log{n}}\right)^{\varepsilon^8\log{n}}\right)\cdot(1-p)^{d_H(s)}\notag\\
&\leq\left(4+(6+o(1))\left(\frac{2e}{\varepsilon^9}\right)^{\varepsilon^8\log{n}}\right)\cdot(1-p)^{d_H(s)}\notag\\
&\leq\left(4+(6+o(1))n^{\varepsilon^8\log{(2e/\varepsilon^9)}}\right)\cdot(1-p)^{d_H(s)}\notag\\
&\leq n^{0.01}\cdot(1-p)^{d_H(s)},
\end{align}
where the second inequality comes from the facts that $p\leq\frac{2\log{n}}{\varepsilon n^{r-1}}$ and $d_H(s)\leq n^{r-1}$ and the last inequality holds since $\varepsilon^8\log{(2e/\varepsilon^9)}$ is sufficiently small whenever $\varepsilon$ is small. 
By symmetry, we deduce
\begin{align}\label{eqn8}
\mathbb{P}[t\in \mathrm{SMALL}(H_{m_3})\,|\,E(T_i)\subseteq E(H_{m_3}), \Sigma(e(s, t), k)]\leq n^{0.01}\cdot(1-p)^{d_H(t)}.
\end{align}

On the other hand, the probability of the event $\{E(T_i)\subseteq E(H_{m_3})\}$ is at most
\begin{align}\label{eqn9}
\mathbb{P}[E(T_i)\subseteq E(H_{m_3})]\leq\frac{\binom{N-|E(T_i)|}{m_3-|E(T_i)|}}{\binom{N}{m_3}}\leq\left(\frac{m_3}{N}\right)^{|E(T_i)|}.
\end{align}

By combining \eqref{eqn6}, \eqref{eqn7}, \eqref{eqn8}, and \eqref{eqn9}, we have
\begin{align*}
\mathbb{P}[s, t\in \mathrm{SMALL}(H_{m_3}), E(T_i) \subseteq E(H_{m_3})]&=\sum_{k=0}^{10r}\mathbb{P}[s\in \mathrm{SMALL}(H_{m_3})\,|\,E(T_i) \subseteq E(H_{m_3}), \Sigma(e(s, t), k)]\\
&\cdot\mathbb{P}[t\in \mathrm{SMALL}(H_{m_3})\,|\,E(T_i) \subseteq E(H_{m_3}), \Sigma(e(s, t), k)]\\
&\cdot\mathbb{P}[\Sigma(e(s, t), k)\,|\,E(T_i) \subseteq E(H_{m_3})]\cdot\mathbb{P}[E(T_i) \subseteq E(H_{m_3})]\\
&+o(n^{-10})\cdot\mathbb{P}[E(T_i) \subseteq E(H_{m_3})]\\
&\leq\left(\sum_{k=0}^{10r}\mathbb{P}[\Sigma(e(s, t), k)\,|\,\Sigma(E(T_i))]\right)\cdot n^{0.02}(1-p)^{d_H(s)+d_H(t)}\cdot p_3^{|E(T_i)|}\\
&+o(n^{-10})\cdot p_3^{|E(T_i)|}\\
&\leq\left( n^{0.02}(1-p)^{d_H(s)+d_H(t)}+o(n^{-10})\right)\cdot p_3^{|E(T_i)|}.
\end{align*}
We now take the union bound over all choices of $T_i$. 
We first sum over all choices of $s, t\in V(H)$, and then sum over all choices of the remaining vertices and edges in $T_i$. By summing over all choices of $s, t\in V(H)$, we have 
\begin{align*}
\sum_{s, t \in V(H)}\mathbb{P}[s, t\in \mathrm{SMALL}(H_{m_3})\,|\,E(T_i) \subseteq E(H_{m_3})] \leq n^{0.02}\left(\sum_{v\in V(H)}(1-p)^{d_H(v)}\right)^2\cdot p_3^{|E(T_i)|}+o(n^{-8})\cdot p_3^{|E(T_i)|}\tag{$**$}\label{stars}
\end{align*}

Since $p=p_3-\frac{O(1)}{n^{r-1}}$, by applying \ref{L2} with $\gamma' = \gamma - O(1)$, we have
\begin{align*}
\sum_{v\in V(H)}(1-p)^{d_H(v)}\leq n^{0.01}.
\end{align*}
Thus, 
\begin{align*}
\eqref{stars}\leq n^{0.05}\cdot p_3^{|E(T_i)|}.
\end{align*}

For each case $i=2, 3, 6, 10$, the number of ways to choose $T_i$, except the choice of $s, t,$ is at most $n^{|E(T_i)|(r-1)-1}$. Thus, by the union bound, the probability that case $i$ occurs is at most $n^{|E(T_i)|(r-1)-1} \cdot n^{0.05}\cdot p_3^{|E(T_i)|} \leq o(1)$. \newline

It now suffices to show that after adding $m_4-m_3$ edges uniformly at random to $H_{m_3}$, $\mathbf{(P3-H_{m_3})}$ still holds. Conditioning on \ref{(P1)} and \ref{(P2)}, the number of vertices that have distance at most $4$ from $\mathrm{SMALL}(H_{m_3})$ is at most $n^{0.1}\cdot\left(\frac{10r}{\varepsilon}\log{n}\right)^4<n^{0.2}$. The probability that at least one edge in this set among the $m_4-m_3$ added edges is at most 
\begin{align*}
(m_4-m_3)\frac{\binom{n^{0.2}}{2} \cdot n^{r-2}}{0.9N}=o(1),
\end{align*}
by the union bound. Thus, after adding the $m_4-m_3$ edges, $\mathbf{(P3-H_{m_3})}$ still holds with high probability.

    \item{\ref{(P4)}} 
For a fixed $U$, there are at most $n^{r-2}|U|^2$ edges in $H$ that intersect $U$ at more than one vertex. Thus, we have 
\begin{align*}
\mathbb{P}[H_{p_2}\notin \mathbf{(P4)}] & \leq \sum_{u=1}^{n/\sqrt{\log{n}}}\binom{n}{u}\cdot \mathbb{P}[\mathrm{Bin}(n^{r-2}u^2, p_2)\geq u(\log{n})^{3/4}]\\
& \stackrel{\eqref{chernoff-U1}}{\leq}\sum_{u=1}^{n/\sqrt{\log{n}}}\binom{n}{u}\cdot \exp\left(-\left(\frac{2u^2\log{n}}{\varepsilon n}\right)\cdot\left(\frac{\varepsilon n}{2u(\log{n})^{1/4}}\right)\log{\left(\frac{\varepsilon n}{6u(\log{n})^{1/4}}\right)}\right)\\
& \leq\sum_{u=1}^{n/\sqrt{\log{n}}}\left(\frac{en}{u}\right)^u\cdot \exp\left(-u(\log{n})^{3/4}\log{\left(\frac{\varepsilon n}{6u(\log{n})^{1/4}}\right)}\right)\\
& \leq\sum_{u=1}^{n/\sqrt{\log{n}}}\exp\left(u\log{(en/u)}-u(\log{n})^{3/4}\log{\left(\frac{\varepsilon n}{6u(\log{n})^{1/4}}\right)}\right)\\
& \leq\sum_{u=1}^{n/\sqrt{\log{n}}}\exp\left(u(1+\log(n/u))-u(\log{n})^{3/4}\log{n}+u(\log{n})^{3/4}\log{\left(\frac{6u(\log{n})^{1/4}}{\varepsilon}\right)}\right)\\
& \leq\sum_{u=1}^{n/\sqrt{\log{n}}}\exp\left(u(\log{n})^{3/4}\left(-0.9\log{n}+\log{\left(\frac{6(\log{n})^{1/4}}{\varepsilon}\right)}\right)\right)\\
& \leq\sum_{u=1}^{n/\sqrt{\log{n}}}\exp\left(-0.1u(\log{n})^{3/4}\cdot\log{n}\right)\\
& \leq2\cdot \exp\left(-0.1(\log{n})^{7/4}\right)=o(1),
\end{align*}
where the last inequality holds by the sum of a geometric series.
    \item{\ref{(P5)}}
For a fixed $U, W$, there are at most $n^{r-2}|U||W|$ edges in $H$ that intersect $U$ at exactly one vertex and also intersect $W$. Without loss of generality, we may assume $|W|=3|U|$. Then the probability that $H_{p_2}\notin \mathbf{(P5)}$ can be bounded by
\begin{align*}
\mathbb{P}[H_{p_2}\notin \mathbf{(P5)}] & \leq \sum_{u=1}^{\frac{n}{\sqrt{\log{n}}}}\binom{n}{u}\binom{n}{3u}\cdot \mathbb{P}[\mathrm{Bin}(3n^{r-2}u^2, p_2)\geq 0.5\varepsilon u\log{n}]\\
& \stackrel{\eqref{chernoff-U1}}{\leq}\sum_{u=1}^{\frac{n}{\sqrt{\log{n}}}}\binom{n}{u}\binom{n}{3u}\cdot \exp\left(-\left(\frac{6u^2\log{n}}{\varepsilon n}\right)\cdot\left(\frac{\varepsilon^2 n}{12u}\right)\log{\left(\frac{\varepsilon^2 n}{36u}\right)}\right)\\
& \leq\sum_{u=1}^{\frac{n}{\sqrt{\log{n}}}}\left(\frac{en}{u}\right)^u\left(\frac{en}{3u}\right)^{3u}\cdot \exp\left(-0.5\varepsilon u\log{n}\cdot\log{\left(\frac{\varepsilon^2 n}{36u}\right)}\right)\\
& \leq\sum_{u=1}^{\frac{n}{\sqrt{\log{n}}}}\exp\left(u\log{\left(\frac{en}{u}\right)}+3u\log{\left(\frac{en}{3u}\right)}-0.5\varepsilon u\log{n}\cdot\log{\left(\frac{\varepsilon^2 n}{36u}\right)}\right)\\
& \leq\sum_{u=1}^{\frac{n}{\sqrt{\log{n}}}}\exp\left(-u\cdot 0.01\varepsilon\log{n}\log{\log{n}}\right)\\
& \leq 2\cdot \exp(-0.01\varepsilon\log{n}\log{\log{n}})=o(1),
\end{align*}
where the last inequality holds by the sum of a geometric series.

    \item{\ref{(P6)}} 
For a fixed $U, W$, and any $u\in U$, the number of edges in $H$ that are incident to $u$ and intersect $W$ at fewer than $r-1$ vertices is at most $\frac{\varepsilon}{2}n^{r-1}$ in $H$. Since the degree of $u$ in $H$ is at least $\left(\frac{1}{2}+\varepsilon\right)\varepsilon n^{r-1}$, the number of edges in $H$ that contain $u$ and intersect $W$ at exactly $r-1$ vertices is at least $\varepsilon^2 n^{r-1}$ in $H$. Thus, there are at least $\varepsilon^2 n^{r-1}|U|=\frac{\varepsilon^2 n^r}{\sqrt{\log{n}}}$ edges that intersect $U$ at exactly one vertex and $W$ at exactly $r-1$ vertices. Therefore, the probability that $H_{p_1}\notin \mathbf{(P6)}$ can be bounded by
\begin{align*}
\mathbb{P}[H_{p_1}\notin \mathbf{(P6)}] & \leq \binom{n}{\frac{n}{\sqrt{\log{n}}}}\binom{n}{(1-\frac{\varepsilon}{2})n}\mathbb{P}[\mathrm{Bin}\left(\frac{\varepsilon^2 n^r}{\sqrt{\log{n}}}, p_1\right)\leq n(\log{n})^{1/3}]\\
& \stackrel{\eqref{chernoff-L1}}{\leq}\binom{n}{\frac{n}{\sqrt{\log{n}}}}\binom{n}{(1-\frac{\varepsilon}{2})n}\cdot \exp\left(-\frac{\varepsilon^2}{3} n\sqrt{\log{n}}\right)\\
& \leq(e\sqrt{\log{n}})^{n/\sqrt{\log{n}}}\left(\frac{e}{1-\frac{\varepsilon}{2}}\right)^{(1-\frac{\varepsilon}{2})n}\cdot \exp\left(-\frac{\varepsilon^2}{3} n\sqrt{\log{n}}\right)=o(1).
\end{align*}

    \item{\ref{(P7)}}
For fixed $U$, $W$, and $u\in U$, since $|U|\leq \frac{n}{2}$, there are at least $\varepsilon n$ vertices $w\in W$ such that $d_H(u, w)>\varepsilon n^{r-2}$. 
Thus, there are at least $\varepsilon^2 n^{r-1}$ edges in $H$ that contain $u$ and intersect $W$. Therefore, there are at least $\frac{\varepsilon^3}{6r}n^{r}$ edges that intersect both $U$ and $W$ in $H$ (each edge can be counted at most $r-1$ times). By taking the union bound over all choices of $U$, we have
\begin{align*}
\mathbb{P}[H_{p_1}\notin \mathbf{(P7)}] & \leq 2^n\cdot \mathbb{P}[\mathrm{Bin}\left(\frac{\varepsilon^3}{6r} n^{r}, p_1\right)=0]\\
& \leq 2^n\cdot \exp\left(-\frac{\varepsilon^3}{6r} n\log{n}\right)=o(1),
\end{align*}
\end{description}
\end{proof}

We now define a random spanning subgraph $\Gamma_0$ of $G=H_{\tau_2}$. 
Let $E_v$ be the set of edges incident to $v$ in $G$ and let $F_v\subseteq E_v$ be a random variable defined as follows:  If $v\in \mathrm{SMALL}(G)$, then $F_v=E_v$. If $v\notin \mathrm{SMALL}(G)$, $F_v$ is a subset of $E_v$ of size exactly $\varepsilon^8\log{n}$, chosen uniformly at random from all such subsets. We set the edge set of $\Gamma_0$ to be $E(\Gamma_0)=\bigcup_{v\in V(G)}F_v$. 

\begin{lemma}\label{lem: expander}
With high probability, $\Gamma_0$ is an $(\frac{\varepsilon n}{6}, 2)$-expander with at most $\varepsilon^8 n\log{n}$ edges.    
\end{lemma}

\begin{proof}
By the definition of $\Gamma_0$, it has at most $\varepsilon^8 n\log{n}$ edges deterministically. Thus, it suffices to check whether $\Gamma_0$ is an $(\varepsilon n/6, 2)$-expander. We assume that~\ref{(P1)}-\ref{(P7)} hold for $G$.
Let $U$ be a subset of $V(\Gamma_0)$ and let $Y$ be a set disjoint from $U$ that covers $U$ (that is, every edge with exactly one endpoint in $U$ also has an endpoint in $Y$). 
We will show that $|Y|\geq 2|U|$ holds by considering two cases depending on the size of $U$.

\begin{enumerate}
    \item{$\frac{n}{\sqrt{\log{n}}}\leq|U|\leq\frac{\varepsilon n}{6}$.}

Let $W=V(G)-(U\cup Y)$. Assume $|Y|< 2|U|$. Then $|W|\geq\left(1-\frac{\varepsilon}{2}\right)n$. For a vertex $u\in U$, let $W_u\subseteq E_u$ be the set of edges which contain $u$ and intersect $W$ in exactly $r-1$ vertices. If there exists $u \in U$ with $|E_u|\leq\varepsilon^8\log{n}$ and $W_u \neq \emptyset$, then $F_u\cap W_u\neq \emptyset$ deterministically. Thus, this contradicts the fact that $Y$ covers $U$.
Let $u$ be a vertex with $|E_u|\geq\varepsilon^8\log{n}$. Note that   by~\ref{(P1)}, $|E_u|\leq\frac{10}{\varepsilon}\log{n}$. Thus, the probability that $F_u\cap W_u=\emptyset$ is at most 
\begin{align*}
\mathbb{P}[F_u\cap W_u=\emptyset]&\leq \frac{\binom{|E_u|-|W_u|}{\varepsilon^8\log{n}}}{\binom{|E_u|}{\varepsilon^8\log{n}}}\leq\left(\frac{|E_u|-|W_u|}{|E_u|}\right)^{\varepsilon^8\log{n}}\leq\exp\left(-10\varepsilon^7(\log{n})^2|W_u|\right).\end{align*}

Note that for any distinct $u, v\in U$, $W_{u}$ and $W_{v}$ are disjoint.
Thus, the events $\{F_u\cap W_u=\emptyset\}$ and $\{F_v\cap W_v=\emptyset\}$ are independent for any distinct $u, v\in U$. 
In addition, by~\ref{(P6)}, $\sum_{u \in U} |W_u| \geq n (\log n)^{1/3}$.
Therefore, the probability that $Y$ does not cover $U$ is at most
\begin{align*}
\mathbb{P}[\bigcup_{u \in U}(F_u \cap W_u) = \emptyset] & =\prod_{u\in U}\mathbb{P}[F_u\cap W_u=\emptyset]\leq\exp(-10\varepsilon^7(\log{n})^2\sum_{u\in U}|W_u|)\\
&\leq\exp(-10\varepsilon^7(\log{n})^{7/3}n)=\exp(-\omega(n)).
\end{align*}
There are at most $2^n$ ways to choose $U$ and $Y$.
Therefore, by the union bound, the probability that there exists a set $U$ with $\frac{n}{\sqrt{\log{n}}}\leq|U|\leq\frac{\varepsilon n}{6}$ that is covered by a set $Y$ with $|Y| \leq 2|U|$ is at most
\begin{align*}
4^n\cdot\exp(-\omega(n))=o(1).
\end{align*}

    \item{$|U|\leq\frac{n}{\sqrt{\log{n}}}$.}

We note that~\ref{(P4)} and \ref{(P5)} hold for $\Gamma_0$, because each is a monotonically decreasing property. 
We denote $U_1:=U\cap \mathrm{SMALL}(G)$ and $U_2:=U\setminus U_1$. Let $Y_1$ be the set of vertices in $Y$ adjacent to $U_1$, and let $Y_2:=Y\setminus Y_1$.

Every edge intersecting $U_1$ in exactly one vertex is a subset of $U\cup Y_1$, thus $U_2\cup Y_1$ covers $U_1$. Since $\delta(H_{\tau_2})=2$, and by~\ref{(P3)} for $\Gamma_0$, $2|U_1|\leq|U_2|+|Y_1|$ holds. Now, by~\ref{(P4)} for $\Gamma_0$, there are at least $|U_2|(\varepsilon^8\log{n}-(\log{n})^{3/4})$ edges intersecting $U_2$ in exactly one vertex. Again, by~\ref{(P3)}, there are at least $|U_2|(\varepsilon^8\log{n}-(\log{n})^{3/4}-1)$ edges intersecting $U_2$ in exactly one vertex and not intersecting $U_1\cup Y_1$. Assume $|Y_2|\leq 3|U_2|$. Then by~\ref{(P5)} for $\Gamma_0$, there are at most $0.5\varepsilon\log{n}|U_2|$ edges intersecting $U_2$ in exactly one vertex and also intersecting $Y_2$. On the other hand, 
\begin{align*}
|U_2|(\varepsilon^8\log{n}-(\log{n})^{3/4}-1)<0.5\varepsilon\log{n}|U_2|
\end{align*}
as $\varepsilon$ is sufficiently small. It is a contradiction. Thus, $|Y_2|\geq 3|U_2|$. Therefore, we have
\begin{align*}
|Y|=|Y_1|+|Y_2|\geq2|U_1|-|U_2|+3|U_2|=2|U|.
\end{align*}
\end{enumerate}    
\end{proof}

We now prove \Cref{lem31}.
\begin{proof}[Proof of \Cref{lem31}]
By \Cref{lem: expander}, with high probability, $G$ contains a spanning subgraph $\Gamma_0$ which is an $(\frac{\varepsilon n}{6}, 2)$-expander with at most $2\varepsilon^8 n\log{n}$ edges.
If $\Gamma_0$ is connected, then we are done.
If $\Gamma_0$ is not connected, as $\Gamma_0$ is an $(\frac{\varepsilon n}{6}, 2)$-expander, the size of each component of $\Gamma_0$ is at least $\frac{\varepsilon n}{6}$. 
Thus, $\Gamma_0$ has at most $\frac{6}{\varepsilon}$ components. 
Let $C$ be a component of $\Gamma_0$.
Then by~\ref{(P7)}, there is at least one edge intersecting both $C$ and $V(\Gamma_0)\setminus C$. Let $\Gamma_1$ be the graph obtained from $\Gamma_0$ by adding that edge. Then the number of components of $\Gamma_1$ is strictly smaller than the number of components of $\Gamma_0$. We repeat this operation at most $\frac{6}{\varepsilon}$ times, and obtain a connected graph.
\end{proof}

\subsection{Proof of \Cref{lem32}}\label{sec: proof of lem32}

We now prove \Cref{lem32}. 
We first prove that any $(\frac{\varepsilon n}{6}, 2)$-expander $\Gamma$ with at most $2\varepsilon^8n\log{n}$ edges has many booster pairs in $H$ (\Cref{lem: many BP}). 
After proving this, we use the fact that there are sufficiently many BPs of $\Gamma$ to prove that $G=H_{\tau_2}$ also has at least one BP for all expanders $\Gamma$ with few edges w.h.p. (\Cref{lem: exist BP2}).

\begin{lemma}\label{lem: many BP}
Let $\Gamma\subseteq H$ be a non-Hamiltonian $(\frac{\varepsilon n}{6}, 2)$-expander with at most $2\varepsilon^8n\log{n}$ edges. Then there is a set $\mathcal{B}$ of booster pairs of $\Gamma$, satisfying the following properties:
\begin{enumerate}
    \item{}
For any $e\in E(H)$, there are at most $\frac{1}{6}\varepsilon^3n^{r-1}$ BPs in $\mathcal{B}$ that contain $e$.
    \item{}
    $|\mathcal{B}|\geq\frac{1}{216r^5}\varepsilon^7n^{2r-1}$.
\end{enumerate}
\end{lemma}

\begin{proof}
Let $P$ be a longest Berge path in $\Gamma$, and denote  $P$ by $(s=p_0, e_0, p_1, \dots, e_{\ell-1}, p_{\ell}=t)$. Let $A_s$ be the set of vertices $v\in V(H)$ such that $d_H(s, v)>\varepsilon n^{r-2}$, and define $A_t$ similarly. By the assumption, $|A_s|, |A_t|\geq\left(\frac{1}{2}+\varepsilon\right)n$. 
Let $B_s\subseteq A_s$ be the set of vertices not contained in $P$, and similarly define $B_t$. Depending on the sizes of $B_s$ and $B_t$, we consider two cases.
\begin{description}
    \item{\textbf{Case 1.}$\enspace|B_s|\geq\frac{2}{3}\varepsilon n$ or $|B_t|\geq\frac{2}{3}\varepsilon n$.}

Without loss of generality, we may assume $|B_s|\geq\frac{2}{3}\varepsilon n$. 
Then there are at least $\frac{2}{3r}\varepsilon^2n^{r-1}$ edges incident to $s$ that contain a  vertex outside of $P$. 
Note that there are at most $2\varepsilon^8n\log{n}$ edges in $\Gamma$, and there are at most $n^{r-2}$ edges incident to both $s$ and $t$. 
After removing edges included in $\Gamma$ or incident to both $s$ and $t$, at least $\frac{1}{2r}\varepsilon^2n^{r-1}$ edges are left. We pick exactly $\frac{1}{2r}\varepsilon^2n^{r-1}$ of these remaining edges and let such collection be $E_s$. Note that adding any edge from $E_s$ to $\Gamma$ increases the maximum length of a Berge path in $\Gamma$. Pick exactly $|E_s|$ edges incident to $t$ and not contained in $\Gamma \cup E_s$, and let the collection of such edges be $E_t$. Note that $E_s$ and $E_t$ are disjoint.

A complete bipartite graph $K_{\alpha, \alpha}$ always has a $\beta$-regular spanning subgraph for every $\beta\leq\alpha$. Thus, one can find a set of booster pairs $\mathcal{B}_{s, t}\subseteq E_s\times E_t$ such that for each edge $e\in E(H)$, the number of BPs in $\mathcal{B}_{s, t}$ containing $e \in E_s \cup E_t$ is exactly $\frac{1}{2r}\varepsilon^2n^{r-2}$. Note that $|\mathcal{B}_{s, t}|=\frac{1}{4r^2}\varepsilon^4n^{2r-3}\geq\frac{1}{6r^3}\varepsilon^5n^{2r-3}$ and for every $(e_1, e_2) \in \mathcal{B}_{s, t}$, $s \in e_1$ and $t \in e_2$.

    \item{\textbf{Case 2.}$\enspace|A_s\setminus B_s|, |A_t\setminus B_t|\geq(\frac{1}{2}+\frac{1}{3}\varepsilon)n$.}

Let $I$ be the set of indices $i$ such that $p_i\in A_t\setminus B_t$, and $p_{i+1}\in A_s\setminus B_s$. By the pigeonhole principle, $|I|\geq\frac{2}{3}\varepsilon{n}$. Let $E'_{s, i}$ be the set of edges incident to both $s$ and $p_{i+1}$, not included in $\Gamma$, and not incident to $t$ or $p_i$. Since $\Gamma$ has at most $2\varepsilon^8 n\log{n}$ edges and there are at most $2n^{r-3}$ edges incident to both $s$ and $p_{i+1}$ and incident to at least one of $t, p_{i}$, we have $|E'_{s, i}|>0.5\varepsilon^2 n^{r-2}$. Let $E_{s, i}$ be an arbitrary subset of $E'_{s, i}$ with size exactly $0.5\varepsilon^2 n^{r-2}$. Also, define $E'_{t, i}$ as the set of edges incident to $t$ and $p_{i}$, not included in $\Gamma$, and not incident to $s$ or $p_{i+1}$. Similarly, let $E_{t, i}$ be an arbitrary subset of $E'_{t, i}$ with size exactly $0.5\varepsilon^2 n^{r-2}$. Then by definition, $E_{s, i}, E_{t, i}$ are disjoint.

For any $e_s\in E_{s, i}, e_t\in E_{t, i}$, $C=P\cup e_s\cup e_t$ contains a Berge cycle that covers all the vertices of $V(P)$. 
Thus, any other edge $e$ of $\Gamma$ that intersects both $V(C)$ and $V(H)\setminus V(C)$ makes the longest path of $\Gamma$ longer by removing $e$ from the Berge cycle if it exists. 
Since $\Gamma$ is connected, there exists such an $e$. Thus $\{e_s, e_t\}$ is a booster pair.

Let $\mathcal{B}_{s, t, i}\subseteq E_{s, i} \times E_{t, i}$ be the set of booster pairs such that for each edge $e\in E(H)$, the number of BPs in $\mathcal{B}_{s, t, i}$ containing $e \in E_{s, i} \cup E_{t, i}$ is exactly $\frac{1}{2r^2}\varepsilon^2n^{r-2}$. Then $|\mathcal{B}_{s, t, i}|=\frac{1}{4r^2}\varepsilon^4n^{2r-4}$. 
Let $\mathcal{B}_{s, t}=\cup_{i\in I}\mathcal{B}_{s, t, i}$. Note that any booster pair can belong to $\mathcal{B}_{s, t, i}$ for at most $r$ indices $i$. Since $|I|\geq \frac{2}{3}\varepsilon n$, $|\mathcal{B}_{s, t}|\geq \frac{1}{6r^3}\varepsilon^5n^{2r-3}$.
Note that for each $e\in E(H)$, there are at most $r$ indices $i$ such that there exists a booster pair $B$ containing $e$ and $B\in\mathcal{B}_{s, t, i}$. 
Thus, the number of booster pairs in $\mathcal{B}_{s, t}$ that contain $e$ is at most $\frac{1}{2r}\varepsilon^2n^{r-2}$.
We also note that for every $(e_1, e_2) \in \mathcal{B}_{s, t}$, $s \in e_1$ and $t \in e_2$.
\end{description}

For a fixed $s$, by \Cref{lem21}, there are at least $\frac{\varepsilon n}{6}$ vertices $t_i$ such that there exists a longest path from $s$ to $t_i$. We pick exactly $\frac{\varepsilon n}{6}$ of them. For each $t_i$, again by \Cref{lem21}, there exist exactly $\frac{\varepsilon n}{6}$ vertices $s_{i, j}$ such that there exists a longest path from $s_{i, j}$ to $t_i$. 

Let
\begin{align*}
\mathcal{B}=\bigcup_{1\leq i, j\leq \frac{\varepsilon n}{6}}\mathcal{B}_{s_{i, j}, t_i}.
\end{align*}
Since each BP is counted at most $r^2$ times in $\mathcal{B}$, the size of $\mathcal{B}$ is $|\mathcal{B}|\geq\frac{1}{r^2}\cdot\frac{\varepsilon^2n^2}{36}\cdot\frac{1}{6r^3}\varepsilon^5n^{2r-3}=\frac{1}{216r^5}\varepsilon^7n^{2r-1}$. For an edge $e\in E(H)$, let $e=\{v_1, \dots, v_r\}$. 
Each $v_i$ can play the role of $s$ or $t$, and by path rotation, there are at most $2r\cdot\frac{\varepsilon n}{6}$ pairs $(s', t')$ such that there exists a BP containing $e$ and contained in $\mathcal{B}_{s', t'}$. Thus, the number of BPs containing $e$ that belong to $\mathcal{B}$ is at most $\frac{1}{6}\varepsilon^3 n^{r-1}$. 
Therefore, $\mathcal{B}$ satisfies both properties.
\end{proof}

We now prove that $G=H_{\tau_2}$ contains at least one booster pair for any $(\frac{\varepsilon n}{6}, 2)$-expander $\Gamma$ with at most $2\varepsilon^8 n\log{n}$ edges w.h.p.
Fix an expander $\Gamma \subseteq G$. 
Let $X=(V, E)$ be an auxiliary graph with $V(X)=E(H)$ and $E(X)$ as the set of booster pairs $B$ with respect to $\Gamma$ defined as in the above lemma.

Let $X^p$ be the probability space of induced subgraphs of $X$, where each vertex is retained with probability $p$ and erased with probability $1-p$, independently of all other vertices.

\begin{lemma}
Let $p_3\leq p\leq p_4$. Then $\mathbb{P}[E(X^p)=\emptyset]\leq 2\cdot \exp\left(-0.4\alpha n^rp\right)$ where $\alpha=\varepsilon^7/\left(216r^5\left(\frac{1}{r!}+\frac{\varepsilon^3}{2}\right)\right)$.
\end{lemma}

\begin{proof}
For a subset $S\subseteq V=V(X)$, let a vertex $x\in V\setminus(S\cup N_{X}(S))$ be $S$-useful if $|N_X(x)\setminus(S\cup N_X(S))|\geq\alpha\cdot n^{r-1}$.

\begin{claim}\label{end}
Let $S\subseteq V$ be a subset such that $|S|\leq\alpha n^r$, $|N_X(S)|\leq2\alpha n^r$. Then the set $A\subseteq V(X)\setminus S$ of all $S$-useful vertices has size at least $\alpha n^r$.
\end{claim}

\begin{proof}
Since $\Delta(X)\leq\frac{1}{6}\varepsilon^3n^{r-1}$, there are at most $\frac{\alpha}{2}\varepsilon^3n^{2r-1}$ edges of $X$ with at least one end in $S\cup N_{X}(S)$. Thus 
\begin{align*}
e_X(V\setminus(S\cup N_X(S)))\geq|E(X)|-\frac{\alpha}{2}\varepsilon^3n^{2r-1}\geq \frac{1}{216r^5}\varepsilon^7n^{2r-1}-\frac{\alpha}{2}\varepsilon^3n^{2r-1}.
\end{align*}
On the other hand,
\begin{align*}
e_X(V\setminus(S\cup N_X(S)))  \leq\Delta(X)\cdot|A|+\alpha n^{r-1}\cdot|V| \leq\frac{1}{6}\varepsilon^3n^{r-1}\cdot|A|+\alpha\binom{n}{r}n^{r-1}\leq\frac{1}{6}\varepsilon^3n^{r-1}\cdot|A|+\alpha\cdot\frac{n^{2r-1}}{r!}.
\end{align*}
Thus $|A|\geq\alpha n^r$.
\end{proof}
Let $v_1,\dots, v_m$ be an arbitrary ordering of $V(X)$. Let $S_0=T_0=\emptyset$, and we do the following operation from $i=1$ to $i=m$. For each $i$, we define $S_i$ and $T_i$ as follows:
\begin{align*}
S_i=
\begin{cases}
S_{i-1}\cup\{v_i\}&~\text{if}~ v_i~\text{is}~T_{i-1}\text{-useful}\\
S_{i-1}&~\text{otherwise}
\end{cases}
\end{align*}
and
\begin{align*}
T_i=
\begin{cases}
T_{i-1}\cup\{v_i\}&~\text{if}~v_i~\text{is}~T_{i-1}\text{-useful and}~v_i\in X^p\\
T_{i-1}&~\text{otherwise}
\end{cases}
\end{align*}

Equivalently, $T_i = S_i \cap V(X^p)$. 
We repeat this operation until $|S_i|\geq\alpha n^r$ or $|N_X(T_i)|\geq2\alpha n^r$.
If $|T_m| \leq |S_m| < \alpha n^r$ and $|N_X(T_m)| < 2\alpha n^r$, then by \Cref{end}, there exist at least $\alpha n^r$ $T_m$-useful vertices.
If a vertex $v_j$ is $T_m$-useful, then it is $T_j$-useful. 
Thus, such a $v_j$ is contained in $S_m$. As $|S_m|<\alpha n^r$, this is a contradiction.
Therefore, the process always ends with either $|S_i|\geq\alpha n^r$ or $|N_X(T_i)|\geq2\alpha n^r$ for some $i \in [m]$. 
Let $t$ be an integer such that the operation terminates at step $t$, and let $S=S_t, T=T_t$. 
Let $\Sigma$ be the event that the operation ends with $|N_X(T)|\geq2\alpha n^r$. 
We note that by the definition of $T_i$-useful, $|N_X(T_i)| \geq |T_i| \cdot \alpha n^{r-1}$ for all $i \in [m]$.
Thus, if $\overline{\Sigma}$ holds, then $|T| < 2n$. 
In addition, we observe that in the $i$-th step, the event that $v_i$ is $T_{i-1}$-useful depends on $X^p$ but conditioned on the event that $v_i$ is $T_{i-1}$-useful, the event that $v_i$ is in $T_i$ is independent of $v_j \in V(X^p)$ for all $j \neq i$.
On the other hand, if $\overline{\Sigma}$ holds, then $|S| \geq \alpha n^r$.
Therefore, the probability of $\overline{\Sigma}$ is at most
\begin{align*}
    \mathbb{P}[\overline{\Sigma}]\leq \mathbb{P}\left[\mathrm{Bin}\left(\alpha n^r, p\right)<2n\right]\stackrel{\eqref{chernoff-L1}}{\leq}\exp\left(-0.4\alpha n^rp\right).
\end{align*}
We now bound the probability that $X^p$ does not contain an edge, conditioned on $\Sigma$. 
We fix a vertex $v_j \notin S$.
Then either $j>t$ or $v_j$ is not $S_{j-1}$-useful. Thus, in the construction of the process, the event that $v_j \in V(X^p)$ does not affect the choice of $T \cap \{v_{j+1}, \ldots, v_m\}$.
Therefore, the event that $v_j \in V(X^p)$ occurs with probability $p$ for any condition on a given $T$. In addition, all such events are mutually independent under the condition on a given $T$. 
If $\Sigma$ holds, then $|S| \leq \alpha n^r$ and $|N_X(T)| \geq 2 \alpha n^r$. Thus, the probability that $N_X(T)\cap V(X^p)=\emptyset$ is
\begin{align*}
\mathbb{P}[N_X(T)\cap V(X^p)=\emptyset\,|\,\Sigma] &\leq \mathbb{P}\left[\mathrm{Bin}\left(\alpha n^r, p\right)=0\right]\stackrel{\eqref{chernoff-L1}}{\leq} \exp\left(-0.4\alpha n^rp\right).
\end{align*}
Note that $E(X^p)=\emptyset$ implies that either $\overline{\Sigma}$ or $\Sigma\wedge\{N_X(T)\cap V(X^p)=\emptyset\}$ holds. Therefore, the probability of $\{E(X^p)=\emptyset\}$ is at most
\begin{align*}
\mathbb{P}[E(X^p)=\emptyset]&\leq\mathbb{P}[\overline{\Sigma}]+\mathbb{P}[N_X(T)\cap V(X^p)=\emptyset\,|\,\Sigma]\cdot\mathbb{P}[\Sigma]\\
&\leq\mathbb{P}[\overline{\Sigma}]+\mathbb{P}[N_X(T)\cap V(X^p)=\emptyset\,|\,\Sigma]\\
&\leq 2\cdot \exp\left(-0.4\alpha n^rp\right).
\end{align*}
\end{proof}

The next step is proving that $H_p$ has a booster pair for not only a specific expander, but also all $(\frac{\varepsilon n}{6}, 2)$-expander with at most $2\varepsilon^8n\log{n}$ edges. 
\begin{lemma}\label{lem: exist BP2}
Let $p_3\leq p\leq p_4$. Then $H_p$ has a booster pair for any $(\frac{\varepsilon n}{6}, 2)$- expander with at most $2\varepsilon^8n\log{n}$ edges with probability at least $1-n^{-\omega(1)}$.
\end{lemma}

\begin{proof}
Let $\mathcal{E}$ be the event that there exists an $(\frac{\varepsilon n}{6}, 2)$-expander $\Gamma$ with at most $2\varepsilon^8\log{n}$ edges such that $H_p$ does not have a booster pair with respect to $\Gamma$. Then $\mathcal{E}$ occurs when there exists $\Gamma$ such that $\Gamma \subseteq H_p$ and $H_p$ has no booster pair with respect to $\Gamma$. By taking the sum over all possible choices of $\Gamma$, we have
\begin{align*}
    \mathbb{P}[\mathcal{E}] & \leq\sum_{i=1}^{2\varepsilon^8n\log{n}}\binom{N}{i}p^i\cdot 2\exp\left(-0.4\alpha n^rp\right)\\
    & \leq\sum_{i=1}^{2\varepsilon^8n\log{n}}\left(\frac{en^rp}{r!i}\right)^i\cdot 2\exp\left(-0.4\alpha n^rp\right).\tag{$***$}\label{starss}
\end{align*}
Since $\left(\frac{en^rp}{r!(i+1)}\right)^{i+1}/\left(\frac{en^rp}{r!i}\right)^{i}=\frac{en^rp}{r!}\cdot\frac{1}{(1+1/i)^i(i+1)}>2$ for $i\leq 2\varepsilon^8n\log{n}$, \eqref{starss} can be bounded by
\begin{align*}
\eqref{starss} & \leq 4\varepsilon^8n\log{n}\left(\frac{en^rp}{2r!\varepsilon^8n\log{n}}\right)^{2\varepsilon^8n\log{n}}\cdot2\exp(-0.4\alpha n^rp)\\
& \leq 8\varepsilon^8 n\log{n}\left(\frac{e}{r!\varepsilon^9}\right)^{2\varepsilon^8 n\log{n}}\cdot \exp(-0.4\alpha n\log{n})\\
& \leq 8\varepsilon^8 n\log{n}\cdot \exp\left((2\varepsilon^8\log{\frac{e}{r!\varepsilon^9}}-0.4\alpha)n\log{n}\right)\\
&=n^{-\omega(1)},
\end{align*}

as $\alpha=\varepsilon^7/\left(216r^5\left(\frac{1}{r!}+\frac{2\varepsilon^3}{3}\right)\right)$ and $\varepsilon$ is small enough.
\end{proof}
We now complete the proof of \Cref{lem32}.
\begin{proof}[Proof of \Cref{lem32}]
Let $\mathcal{E}_m$ be the event that $H_m$ contains an $(\frac{\varepsilon n}{6}, 2)$-expander with at most $2\varepsilon^8n\log{n}$ edges for which $H_m$ does not have a booster pair. Since $m_3\leq\tau_2\leq m_4$ with probability $1-o(1)$ by \Cref{lem26}, 
\begin{align*}
\mathbb{P}[\mathcal{E}_{\tau_2}] \leq o(1)+\sum_{m=m_3}^{m_4}\mathbb{P}[\mathcal{E}_m]\leq o(1)+3(m_4-m_3)\cdot n^{-\omega(1)}=o(1).
\end{align*}
\end{proof}

\section{Proof outline of \Cref{theorem: berge} for weak Berge Hamiltonicity}\label{sec:weak_berge}

In this section, we outline the proof of the case of a weak Berge Hamilton cycle.

The first step is estimating the value of $\tau_1$ as in \Cref{sec:value_of_tau}. 
We define $q_0$ as the unique real root of the following equation: $$\sum_{v\in V(H)}(1-p)^{d_H(v)}=1.$$ Note that the left-hand side of this equation is strictly decreasing as $p$ increases. Let $\gamma=\gamma(n)$ be a nonnegative increasing function such that $\gamma=\Theta(\log{\log{\log{n}}})$. Let 
\begin{align*}
p_5=q_0-\frac{\gamma}{n^{r-1}}, p_6=q_0+\frac{\gamma}{n^{r-1}}, m_5=N\cdot p_5, m_6=N\cdot p_6.
\end{align*}
Then one can prove that $p_5\leq\tau_1\leq p_6$ w.h.p.
Based on this estimation of $\tau_1$, together with  $p_1\leq\tau_1\leq p_2$, we can prove that $H_{\tau_1}$ satisfies the properties~\ref{(P1)}-\ref{(P7)} with high probability. 

The rest of the proof is finding an expander with few edges in $H_{\tau_1}$ and proving that $H_{\tau_1}$ contains a booster pair for any such expander.
On the other hand, the previous definition of a $(k, \alpha)$-expander is not appropriate for the weak Berge problem.
The following definition of an expander and the corresponding lemma is introduced in~\cite{deepak}.
An $r$-graph $H$ is a \emph{weak $(k, \alpha)$-expander} if $|N_H(X)| \geq \alpha |X|$ for any set of vertices $X$ with $|X|\leq k$.
\begin{lemma}\label{lem:rotation_weak}
    Let $r \geq 3$ and $H$ be an $r$-graph on $n$ vertices which is a weak $(k, 2)$-expander. Let $P = (v_1, e_1, \ldots, e_{\ell-1}, v_\ell)$ be a longest weak Berge path in $H$ with one endpoint $v_1=x$.
    Then there are at least $k$ distinct longest weak Berge paths in $H$ that can be obtained from a sequence of path rotations fixing $x$.
\end{lemma}

We now take a spanning subgraph $\Gamma_0$ of $H_{\tau_1}$ by choosing $\varepsilon^8 \log{n}$ random edges incident to each vertex $v \in V(H)$ if $v \notin \mathrm{SMALL}(H_{\tau_1})$, and all edges incident to $v$ if $v \in \mathrm{SMALL}(H_{\tau_1})$.  
To show that $\Gamma_0$ is a weak $(\frac{\varepsilon n}{6}, 2)$-expander with high probability, the proof is almost identical to the proof of \Cref{lem: expander}.
We define $Y$ as $N_{\Gamma_0}(U)$ instead of a set that covers $U$ and claim that $|Y|\geq 2|U|$.
Then, the only difference is in the case where $|U|\leq\frac{n}{\sqrt{\log{n}}}$. In this case, we define $U_1 = U \cap \mathrm{SMALL}(H_{\tau_1})$ and $Y_1 = Y \cap N_{H_{\tau_1}}(U_1)$. 
Since $r \geq 3$ and $\delta(H_{\tau_1})=1$, we have $2|U_1|\leq |N_H(U_1)| \leq |Y_1| + |U_2|$ where $U_2 = U \setminus U_1$. 
The rest of the proof is identical to that of \Cref{lem: expander}.

Finally, by applying \Cref{lem:rotation_weak} instead of \Cref{lem21}, we can find many booster pairs for any weak $(\frac{\varepsilon n}{6}, 2)$-expander with at most $2\varepsilon^8 n\log{n}$ edges. Then, by the same argument as in \Cref{sec: proof of lem32}, we can show that an analogous version of \Cref{lem32} holds for the weak Berge case, which completes the proof of the weak Berge result.

\end{document}